\documentclass[twoside]{amsart}

\newif\ifdraft
\draftfalse
\usepackage{amsmath,amsfonts,amsthm,mathrsfs}
\usepackage{amssymb}

\usepackage[ocgcolorlinks,unicode,bookmarks]{hyperref}
\usepackage[usenames,dvipsnames]{xcolor}
\hypersetup{colorlinks=true,citecolor=NavyBlue,linkcolor=Brown,urlcolor=Green} 

\usepackage{expl3}

\ExplSyntaxOn
\prop_new:N \g_cite_map_prop
\tl_new:N \l_citekey_result_tl

\cs_new:Npn \mapcitekey #1#2 {
  \clist_map_inline:nn {#2}
       {  \prop_gput:Nnn  \g_cite_map_prop  {##1} {#1}   }
}

\cs_new:Npn \getcitekey #1 {
   \prop_get:NoN \g_cite_map_prop{#1}  \l_citekey_result_tl
   \quark_if_no_value:NF \l_citekey_result_tl
       {  \tl_set_eq:NN #1  \l_citekey_result_tl  }
}

\cs_new:Npn \showcitekeymaps {\prop_show:N  \g_cite_map_prop }
\ExplSyntaxOff

\usepackage{etoolbox}
\makeatletter
\patchcmd{\@citex}{\if@filesw}{\getcitekey\@citeb \if@filesw}%
    {\typeout{*** SUCCESS ***}}{\typeout{*** FAIL ***}}
\patchcmd{\nocite}{\if@filesw}{\getcitekey\@citeb \if@filesw}%
    {\typeout{*** SUCCESS ***}}{\typeout{*** FAIL ***}}
\makeatother

\usepackage{enumitem}
\newenvironment{renumerate}{%
	\begin{enumerate}[label=(\roman{*}), ref=(\roman{*})]
}{%
	\end{enumerate}%
}

\newenvironment{aenumerate}{%
	\begin{enumerate}[label=(\alph{*}), ref=(\alph{*})]
}{%
	\end{enumerate}%
}

\usepackage{chngcntr}

\ifdraft
\usepackage[notcite,notref,color]{showkeys}

\definecolor{labelkey}{gray}{0.5}
\fi

\usepackage{tikz}
\usepackage{tikz-cd}

\newcommand{\MHM}{\operatorname{MHM}}

\newcommand{\Dmod}{\mathscr{D}}
\newcommand{\Mmod}{\mathcal{M}}
\newcommand{\Nmod}{\mathcal{N}}

\newcommand{\shT}{\mathscr{T}}

\newcommand{\derR}{\mathbf{R}}
\newcommand{\derL}{\mathbf{L}}

\newcommand{\decal}[1]{\lbrack #1 \rbrack}

\newcommand{\shH}{\mathcal{H}}

\newcommand{\abs}[1]{\lvert #1 \rvert}

\newcommand{\tensor}{\otimes}

\newcommand{\shHom}{\mathcal{H}\hspace{-1pt}\mathit{om}}

\newcommand{\NN}{\mathbb{N}}
\newcommand{\ZZ}{\mathbb{Z}}
\newcommand{\QQ}{\mathbb{Q}}
\newcommand{\RR}{\mathbb{R}}
\newcommand{\CC}{\mathbb{C}}

\newcommand{\PP}{\mathbb{P}}

\newcommand{\menge}[2]{\bigl\{ \thinspace #1 \thinspace\thinspace \big\vert%
\thinspace\thinspace #2 \thinspace \bigr\}}

\DeclareMathOperator{\im}{im}

\DeclareMathOperator{\Supp}{Supp}
\DeclareMathOperator{\codim}{codim}

\DeclareMathOperator{\Sym}{Sym}
\DeclareMathOperator{\gr}{gr}
\DeclareMathOperator{\DR}{DR}

\DeclareMathOperator{\Var}{Var}

\DeclareMathOperator{\Pic}{Pic}

\newcommand{\define}[1]{\emph{#1}}

\newcommand{\shf}[1]{\mathscr{#1}}
\newcommand{\OX}{\shf{O}_X}
\newcommand{\OmX}{\Omega_X}

\newcommand{\shV}{\shf{V}}

\def\overbar#1#2#3{{%
	\setbox0=\hbox{\displaystyle{#1}}%
	\dimen0=\wd0
	\advance\dimen0 by -#2 
	\vbox {\nointerlineskip \moveright #3 \vbox{\hrule height 0.3pt width \dimen0}%
		\nointerlineskip \vskip 1.5pt \box0}%
}}

\newcommand{\fu}{f^{\ast}}
\newcommand{\fl}{f_{\ast}}

\newcommand{\qu}{q^{\ast}}

\newcommand{\tl}{t_{\ast}}

\newcommand{\hl}{h_{\ast}}

\newcommand{\DD}{\mathbb{D}}

\newcommand{\shF}{\shf{F}}
\newcommand{\shG}{\shf{G}}

\newcommand{\shO}{\shf{O}}

\makeatletter
\let\@@seccntformat\@seccntformat
\renewcommand*{\@seccntformat}[1]{%
  \expandafter\ifx\csname @seccntformat@#1\endcsname\relax
    \expandafter\@@seccntformat
  \else
    \expandafter
      \csname @seccntformat@#1\expandafter\endcsname
  \fi
    {#1}%
}
\newcommand*{\@seccntformat@subsection}[1]{%
  \textbf{\csname the#1\endcsname.}
}
\makeatother

\makeatletter
\let\@paragraph\paragraph
\renewcommand*{\paragraph}[1]{%
	\vspace{0.3\baselineskip}%
	\@paragraph{\textit{#1}}%
}
\makeatother

\counterwithin{equation}{subsection}
\counterwithout{subsection}{section}
\counterwithin{figure}{subsection}

\newtheorem{theorem}[equation]{Theorem}
\newtheorem*{theorem*}{Theorem}

\newtheorem*{lemma*}{Lemma}
\newtheorem{corollary}[equation]{Corollary}
\newtheorem*{corollary*}{Corollary}
\newtheorem{proposition}[equation]{Proposition}
\newtheorem*{proposition*}{Proposition}

\newtheorem*{conjecture*}{Conjecture}

\theoremstyle{definition}
\newtheorem{definition}[equation]{Definition}
\newtheorem*{definition*}{Definition}
\newtheorem{remark}[equation]{Remark}
\newtheorem*{remark*}{Remark}

\newtheorem{example}[equation]{Example}
\newtheorem*{example*}{Example}
\newtheorem*{problem*}{Problem}

\theoremstyle{plain}

\newcommand{\theoremref}[1]{\hyperref[#1]{Theorem~\ref*{#1}}}
\newcommand{\lemmaref}[1]{\hyperref[#1]{Lemma~\ref*{#1}}}
\newcommand{\definitionref}[1]{\hyperref[#1]{Definition~\ref*{#1}}}
\newcommand{\propositionref}[1]{\hyperref[#1]{Proposition~\ref*{#1}}}
\newcommand{\conjectureref}[1]{\hyperref[#1]{Conjecture~\ref*{#1}}}
\newcommand{\corollaryref}[1]{\hyperref[#1]{Corollary~\ref*{#1}}}
\newcommand{\exampleref}[1]{\hyperref[#1]{Example~\ref*{#1}}}
\newcommand{\exerciseref}[1]{\hyperref[#1]{Exercise~\ref*{#1}}}

\makeatletter
\let\old@caption\caption
\renewcommand*{\caption}[1]{%
	\setcounter{figure}{\value{equation}}%
	\stepcounter{equation}%
	\old@caption{#1}\relax%
}
\makeatother

\newcounter{intro}

\newtheorem{intro-conjecture}[intro]{Conjecture}
\newtheorem{intro-corollary}[intro]{Corollary}
\newtheorem{intro-theorem}[intro]{Theorem}

\newcommand{\OmA}{\Omega_A}
\newcommand{\OmB}{\Omega_B}
\newcommand{\Ah}{\hat{A}}
\newcommand{\omY}{\omega_Y}

\newcommand{\OmY}{\Omega_Y}

\newcommand{\df}{\mathit{df}}

\newcommand{\parref}[1]{\hyperref[#1]{\S\ref*{#1}}}
\newcommand{\chapref}[1]{\hyperref[#1]{Chapter~\ref*{#1}}}

\makeatletter
\newcommand*\if@single[3]{%
  \setbox0\hbox{${\mathaccent"0362{#1}}^H$}%
  \setbox2\hbox{${\mathaccent"0362{\kern0pt#1}}^H$}%
  \ifdim\ht0=\ht2 #3\else #2\fi
  }
\newcommand*\rel@kern[1]{\kern#1\dimexpr\macc@kerna}
\newcommand*\widebar[1]{\@ifnextchar^{{\wide@bar{#1}{0}}}{\wide@bar{#1}{1}}}
\newcommand*\wide@bar[2]{\if@single{#1}{\wide@bar@{#1}{#2}{1}}{\wide@bar@{#1}{#2}{2}}}
\newcommand*\wide@bar@[3]{%
  \begingroup
  \def\mathaccent##1##2{%
    \if#32 \let\macc@nucleus\first@char \fi
    \setbox\z@\hbox{$\macc@style{\macc@nucleus}_{}$}%
    \setbox\tw@\hbox{$\macc@style{\macc@nucleus}{}_{}$}%
    \dimen@\wd\tw@
    \advance\dimen@-\wd\z@
    \divide\dimen@ 3
    \@tempdima\wd\tw@
    \advance\@tempdima-\scriptspace
    \divide\@tempdima 10
    \advance\dimen@-\@tempdima
    \ifdim\dimen@>\z@ \dimen@0pt\fi
    \rel@kern{0.6}\kern-\dimen@
    \if#31
      \overline{\rel@kern{-0.6}\kern\dimen@\macc@nucleus\rel@kern{0.4}\kern\dimen@}%
      \advance\dimen@0.4\dimexpr\macc@kerna
      \let\final@kern#2%
      \ifdim\dimen@<\z@ \let\final@kern1\fi
      \if\final@kern1 \kern-\dimen@\fi
    \else
      \overline{\rel@kern{-0.6}\kern\dimen@#1}%
    \fi
  }%
  \macc@depth\@ne
  \let\math@bgroup\@empty \let\math@egroup\macc@set@skewchar
  \mathsurround\z@ \frozen@everymath{\mathgroup\macc@group\relax}%
  \macc@set@skewchar\relax
  \let\mathaccentV\macc@nested@a
  \if#31
    \macc@nested@a\relax111{#1}%
  \else
    \def\gobble@till@marker##1\endmarker{}%
    \futurelet\first@char\gobble@till@marker#1\endmarker
    \ifcat\noexpand\first@char A\else
      \def\first@char{}%
    \fi
    \macc@nested@a\relax111{\first@char}%
  \fi
  \endgroup
}
\makeatother

\newcommand{\omX}{\omega_X}

\newcommand{\shA}{\mathcal{A}}
\renewcommand{\DD}{\mathbf{D}}
\newcommand{\AAA}{\mathbf{A}}

\newcommand{\shGb}{\shG_{\bullet}}
\DeclareMathOperator{\HM}{HM}

\renewcommand{\shV}{\mathcal{V}}

\newcommand{\V}{\mathbb{V}}

\newcommand{\cV}{\mathcal{V}}

\mapcitekey{Campana+Peternell:GeometricStability}{CPet}
\mapcitekey{Campana+Paun:Quotients}{CP2}
\mapcitekey{Campana+Paun:OrbifoldOne}{CP1}
\mapcitekey{Kollar:DualizingI}{Kollar1}
\mapcitekey{Kollar:DualizingII}{Kollar3}
\mapcitekey{Kollar:Subadditivity}{Kollar2}
\mapcitekey{Kollar:Shafarevich}{Kollar4}
\mapcitekey{Viehweg:WP1}{Viehweg1}
\mapcitekey{Viehweg:WP2}{Viehweg2}
\mapcitekey{Viehweg:Survey}{Viehweg3}
\mapcitekey{Viehweg:Hilbert1}{Viehweg4}
\mapcitekey{Popa+Wu:WeakPositivity}{PW}
\mapcitekey{Viehweg+Zuo:Isotriviality}{VZ1}
\mapcitekey{Viehweg+Zuo:BaseSpaces}{VZ2}
\mapcitekey{Viehweg+Zuo:Brody}{VZ3}
\mapcitekey{Zuo:Kernels}{Zuo}
\mapcitekey{Popa+Schnell:OneForms}{PS}
\mapcitekey{Popa+Schnell:Hyperbolicity}{PS2}
\mapcitekey{Popa+Schnell:mhmgv}{mhmgv}
\mapcitekey{Saito:KC}{Saito-KC}
\mapcitekey{BDPP:Pseudoeffective}{BDPP}
\mapcitekey{Esnault+Viehweg:VanishingTheorems}{EV}
\mapcitekey{Saito:HodgeModules}{Saito-MHP}
\mapcitekey{Saito:MixedHodgeModules}{Saito-MHM}
\mapcitekey{Saito:b-function}{Saito-B}
\mapcitekey{Saito:Steenbrink}{Saito-Spectrum}
\mapcitekey{Saito:Kollar}{Saito-Kollar}
\mapcitekey{Saito:YPG}{Saito-YPG}
\mapcitekey{Popa:SaitoVanishing}{Popa}
\mapcitekey{Brunebarbe:Variations}{Brunebarbe}
\mapcitekey{Schnell:HolonomicAbelian}{Schnell}
\mapcitekey{Kebekus+Kovacs:FamiliesSurfaces}{KK1}
\mapcitekey{Kebekus+Kovacs:FamiliesCompact}{KK2}
\mapcitekey{Kebekus+Kovacs:FamiliesThreefolds}{KK3}
\mapcitekey{Kovacs:FamiliesSurvey}{Kovacs}
\mapcitekey{Kovacs:BaseNef}{Kovacs2}
\mapcitekey{Kovacs:Hyperbolicity}{Kovacs3}
\mapcitekey{Kebekus:VZSurvey}{Kebekus}
\mapcitekey{Patakfalvi:Hyperbolicity}{Patakfalvi}
\mapcitekey{Kawamata:AdditivityMMP}{Kawamata}
\mapcitekey{Mori:Survey}{Mori}
\mapcitekey{Schnell:sanya}{Schnell1}
\mapcitekey{Schnell:saito-vanishing}{Schnell2}
\mapcitekey{Schnell:weakpos}{Schnell3}
\mapcitekey{Schnell:Neron}{Schnell4}
\mapcitekey{Schnell:ExtendedHodge}{Schnell5}
\mapcitekey{Taji:special}{Taji}
\mapcitekey{Mustata+Popa:HodgeIdeals}{MP}
\mapcitekey{Pareschi+Popa+Schnell:Kaehler}{PPS}
\mapcitekey{Esnault+Viehweg:VanishingTheorems}{EV}
\mapcitekey{Wu:Vanishing}{Wu}
\mapcitekey{Suh:Vanishing}{Suh}
\mapcitekey{Hotta+Takeuchi+Tanisaki:Book}{HTT}
\mapcitekey{Ein+Lazarsfeld:Theta}{EL}
\mapcitekey{Schmid:VHS}{Schmid}
\mapcitekey{Hacon+Kovacs:OneForms}{HK}
\mapcitekey{Luo+Zhang:Threefolds}{LZ}
\mapcitekey{Zhang:FormsAmple}{Zhang}
\mapcitekey{BCHM:MinimalModels}{BCHM}
\mapcitekey{Saito:Kaehler}{Saito-Kaehler}
\mapcitekey{Green+Lazarsfeld:GV1}{GL1}
\mapcitekey{Green+Lazarsfeld:GV2}{GL2}
\mapcitekey{Hacon:GV}{Hacon}
\mapcitekey{Chen+Jiang:Pushforward}{CJ}
\mapcitekey{Wang:TorsionPoints}{Wang}
\mapcitekey{Pareschi:Standard}{Pareschi}
\mapcitekey{Chen+Hacon:abelian}{CH1}
\mapcitekey{Jiang:effective}{Jiang}
\mapcitekey{Simpson:Subspaces}{Simpson}
\mapcitekey{Pareschi+Popa:GV}{PP3}
\mapcitekey{Pareschi+Popa:regI}{PP4}
\mapcitekey{Pareschi+Popa:regIII}{PP5}
\mapcitekey{Lazarsfeld:PositivityII}{Lazarsfeld}
\mapcitekey{Szendroi:CDTSurvey}{Szendroi}
\mapcitekey{BBDJS:VanishingCycles}{BBDJS}
\mapcitekey{Kiem+Li:CategorificationDT}{KL}
\mapcitekey{Schmid+Vilonen:Unitary}{SV}
\mapcitekey{Achar+Kitchen:KoszulDuality}{AK}
\mapcitekey{Cautis+Dodd+Kamnitzer:sl2action}{CDK}
\mapcitekey{Wu:Multiplier}{Wu2}
\mapcitekey{Fujino+Fujisawa+Saito:Semipositivity}{FFS}
\mapcitekey{Budur+Saito:Multiplier}{BS}
\mapcitekey{Dimca+Maisonobe+Saito:Spectrum}{DMS}
\mapcitekey{Maxim+Schuermann:CharacteristicClasses}{MS}
\mapcitekey{Sabbah+Schnell:MHMProject}{SS}
\mapcitekey{Brosnan+Fang+Nie+Pearlstein:Singularities}{BFNP}

\begin{document}

\title[Positivity for Hodge modules and geometric applications]{Positivity for Hodge modules and geometric applications}

\author{Mihnea Popa}
\address{Department of Mathematics, Northwestern University,
2033 Sheridan Road, Evanston, IL 60208, USA} 
\email{mpopa@math.northwestern.edu}

\thanks{The author was partially supported by the NSF grant DMS-1405516 and a Simons Fellowship.}


\begin{abstract}
This is a survey of vanishing and positivity theorems for Hodge modules, and their recent applications 
to birational and complex geometry.
\end{abstract}


\maketitle

\markboth{MIHNEA POPA} 
{POSITIVITY FOR HODGE MODULES AND GEOMETRIC APPLICATIONS}

\makeatletter
\newcommand\@dotsep{4.5}
\def\@tocline#1#2#3#4#5#6#7{\relax
  \ifnum #1>\c@tocdepth 
  \else
    \par \addpenalty\@secpenalty\addvspace{#2}%
    \begingroup \hyphenpenalty\@M
    \@ifempty{#4}{%
      \@tempdima\csname r@tocindent\number#1\endcsname\relax
    }{%
      \@tempdima#4\relax
    }%
    \parindent\z@ \leftskip#3\relax
    \advance\leftskip\@tempdima\relax
    \rightskip\@pnumwidth plus1em \parfillskip-\@pnumwidth
    #5\leavevmode\hskip-\@tempdima #6\relax
    \leaders\hbox{$\m@th
      \mkern \@dotsep mu\hbox{.}\mkern \@dotsep mu$}\hfill
    \hbox to\@pnumwidth{\@tocpagenum{#7}}\par
    \nobreak
    \endgroup
  \fi}
\def\l@section{\@tocline{1}{0pt}{1pc}{}{\bfseries}}
\def\l@subsection{\@tocline{2}{0pt}{25pt}{5pc}{}}
\makeatother

\tableofcontents


\setlength{\parskip}{0.5\baselineskip}

\section{Introduction}

Over the past few years a substantial body of work has been devoted to applications of Morihiko Saito's theory of mixed Hodge modules to the geometry of complex algebraic varieties, as well as to questions in singularities, arithmetic geometry, representation theory, and other areas.  
One of the main reasons for this is the realization that generalizations of classical constructions in Hodge theory, like Hodge bundles and de Rham complexes, arising in this context, satisfy analogues of well-known vanishing and positivity theorems.

This survey has two main aims. On one hand, I will review the vanishing, injectivity and weak positivity results that are known in the context of Hodge modules; it turns out that there is by now a complete package of such theorems extending the standard results of this type in birational geometry. On the other hand, I will describe applications of these results in, loosely speaking, the realm of birational geometry: results on families of varieties, holomorphic forms, generic vanishing theorems, and singularities, 
that at the moment we only know how to approach using the theory of Hodge modules.

Here is a brief review of the contents of the paper; each section contains a more detailed description of the topic treated there. 
First, in \S\ref{Dmod} and \S\ref{Hmod} I include a quick tour of the basic facts from the theory of filtered $\Dmod$-modules and 
Hodge modules, respectively, that are used in the paper. There are also references to general material on these subjects.

Ch.\ref{vanishing_positivity} reviews the package of vanishing and positivity theorems for Hodge modules mentioned above. 
Ch.\ref{generic_vanishing} contains applications of Hodge modules to generic vanishing. This is an area where their impact has been 
substantial, both by allowing to extend the range of applications in the projective setting, and by providing basic new results in 
the K\"ahler setting.

Both standard and generic vanishing for Hodge modules, as well as weak positivity, have proven to be useful in the 
study of families of varieties. The applications are obtained by means of a general construction, inspired by objects introduced by Viehweg and Zuo in a somewhat more restrictive context, which is described in \S\ref{HM-construction}. In \S\ref{one_forms} I explain how this 
construction is used in studying zeros of holomorphic one-forms, while in \S\ref{abvar} this is applied to families of varieties of general type 
over abelian varieties. On the other hand, \S\ref{VZ} contains applications of this construction to producing Viehweg-Zuo sheaves 
for families with maximal variation whose geometric generic fiber has a good minimal model; this in particular leads to a strengthening
of Viehweg's hyperbolicity conjecture.

Ch.\ref{Hodge_ideals} is devoted to an introduction to recent joint work with M. Musta\c t\u a, whose focus are objects called Hodge ideals,
intimately linked to Saito's Hodge filtration on the $\Dmod$-module of functions with poles along a given divisor. They are 
a useful tool in studying the singularities and Hodge theory of hypersurfaces and their complements in smooth complex varieties.
Both global methods like vanishing theorems, and local techniques from the study of Hodge modules, are important here.

The paper is certainly not intended to be a balanced survey of all recent applications of Hodge modules, which is something beyond my 
abilities. I have however included in Ch.\ref{other_applications} a (non-exhaustive) list, with brief explanations, of notable uses of Hodge modules to areas other than those treated here, and references to further literature.


\noindent
{\bf Acknowledgements.}
This survey is a (much) expanded version of the lecture I gave at the Utah AMS Summer Institute in July 2015. I would like to thank the organizers of the respective seminar, Stefan Kebekus and Claire Voisin, for the invitation and for suggesting that I submit a paper. Many results described here were obtained together with Christian Schnell, whom I would like to thank for numerous discussions and explanations on the theory of Hodge modules. I would also like to thank Mircea Musta\c t\u a  and Lei Wu for useful conversations on the topic.

\section{Preliminaries on Hodge modules}

The next two sections contain the minimal amount of information necessary in order to navigate the rest of the paper, for the reader who is not closely acquainted to the theory of Hodge modules; those who are can skip directly to Ch.\ref{vanishing_positivity}.
The original, and most comprehensive, references for Hodge modules are Saito's papers \cite{Saito-MHP} and \cite{Saito-MHM}. Recently Saito also wrote a more informal introduction in \cite{Saito-YPG}. Another general overview of the theory can be 
found in \cite{Schnell1}. A detailed source in its development stages, but already containing plenty of useful information, is 
the project \cite{SS} by Sabbah and Schnell.

\subsection{Background on $\Dmod$-modules}\label{Dmod}
This section contains a very quick review of the necessary tools from the theory of $\Dmod$-modules. A comprehensive 
reference for the material here is \cite{HTT}.

\noindent
{\bf Filtered $\Dmod$-modules.}
Let $X$ be a  smooth complex variety of dimension $n$. A filtered left $\Dmod$-module on $X$ is a $\Dmod_X$-module $\Mmod$ with an increasing filtration 
$F = F_\bullet \Mmod$ by coherent $\shO_X$-modules, bounded from below and satisfying
$$F_k \Dmod_X \cdot F_\ell \Mmod  \subseteq F_{k+ \ell} \Mmod \,\,\,\,{\rm for~all~} k, \ell \in \ZZ.$$
In addition, the filtration is \emph{good} if the inclusions above are equalities for $k \gg 0$. This condition is equivalent to the fact that the total associated graded object
$$\gr^F_{\bullet} \Mmod = \bigoplus_k \gr_k^F \Mmod = \bigoplus_k F_k \Mmod / F_{k-1} \Mmod$$
is finitely generated over $\gr_{\bullet}^F \Dmod_X \simeq {\rm Sym}~T_X$, i.e. induces a coherent sheaf on the cotangent bundle 
$T^*X$. Similar definitions apply to right $\Dmod_X$-modules. 

While the simplest example of a filtered left $\Dmod_X$-module is $\shO_X$,
with the trivial filtration $F_k \shO_X = \shO_X$ for $k \ge 0$ and $F_k \shO_X = 0$ for $k < 0$, the canonical bundle $\omega_X$ is naturally endowed with a filtered right $\Dmod_X$-module structure; locally, if $z_1, \ldots, z_n$ are local coordinates on $X$, for any $f \in \shO_X$ and any $P \in \Dmod_X$ the action is
$$(f \cdot dz_1 \wedge \cdots \wedge dz_n) \cdot P = {}^tP (f) \cdot dz_1 \wedge \cdots \wedge dz_n.$$ 
Here, if $P = \sum_{\alpha} g_\alpha \partial^\alpha$, then ${}^t P =  \sum_{\alpha} (-\partial)^\alpha g_\alpha$ is its formal adjoint.
This structure leads to a one-to-one correspondence between left and right $\Dmod_X$-modules given by 
$$\Mmod \mapsto \Nmod = \Mmod \otimes{_{\shO_X}} \omega_X \,\,\,\, {\rm and} \,\,\,\, \Nmod \mapsto \Mmod = \mathcal{H}om_{\shO_X} (\omega_X, \Nmod),$$
while on filtrations the left-right rule is
$$F_p \Mmod = F_{p- n} \Nmod \otimes_{\shO_X} \omega_X^{-1}.$$
In particular $F_k \omega_X = \omega_X$ for $k \ge - n$ and $0$ for $k < - n$.

\noindent
{\bf de Rham complex.}
Given a left $\Dmod_X$-module $\Mmod$, its de Rham complex is:
\[
\DR (\Mmod) = \Bigl\lbrack
		\Mmod \to \Omega_X^1 \tensor \Mmod \to \dotsb \to \Omega_X^n \tensor \Mmod
	\Bigr\rbrack.
\]
This is a $\CC$-linear complex with differentials induced by the corresponding integrable connection
$\nabla: \Mmod \rightarrow \Mmod \otimes \Omega_X^1$. By placing it in degrees $-n, \ldots, 0$,
we consider it to be the de Rham complex associated to the corresponding right $\Dmod$-module.
A filtration $F_{\bullet} \Mmod$ on $\Mmod$ induces a  filtration on $\DR (\Mmod)$ by the formula
$$
	F_k \DR(\Mmod) = \Bigl\lbrack
		F_k \Mmod \to \Omega_X^1 \otimes F_{k+1} \Mmod \to \dotsb 
			\to \Omega_X^n \otimes F_{k+n} \Mmod
	\Bigr\rbrack.
$$
For any integer $k$, the associated graded complex for this filtration is 
$$
	\gr_k^F \DR(\Mmod) = \Big\lbrack
		\gr_k^F \Mmod \to \Omega_X^1 \otimes \gr_{k+1}^F \Mmod \to \dotsb \to
			\Omega_X^n \otimes \gr_{k+n}^F \Mmod
	\Big\rbrack.
$$
This is now a complex of coherent $\OX$-modules in degrees $-n, \dotsc, 0$, 
providing an object in ${\bf D}^b (X)$, the bounded derived category of coherent sheaves on $X$.

The lowest non-zero graded piece of a filtered $\Dmod$-module is of particular interest.
For one such left $\Dmod_X$-module $(\Mmod, F)$ define
$$
p (\Mmod) : = {\rm min}~\{p ~|~ F_p \Mmod \neq 0\}\,\,\,\,{\rm and} \,\,\,\, S (\Mmod) := F_{p(\Mmod)} \Mmod.
$$
For the associated right $\Dmod_X$-module we then have
$$p(\Nmod) = p (\Mmod) - n \,\,\,\, {\rm and} \,\,\,\,S(\Nmod) = S(\Mmod) \otimes \omega_X. $$

\noindent
{\bf Pushforward.}
Let $f\colon X \rightarrow Y$ be a morphism of smooth complex varieties. The associated transfer module
$$\Dmod_{X\to Y} : = \shO_X \otimes_{f^{-1} \shO_Y} f^{-1} \Dmod_Y$$
has the structure of a $(\Dmod_X, f^{-1} \Dmod_Y)$-bimodule, and has a filtration given by 
$f^* F_k \Dmod_Y$. The pushforward via $f$ is a natural operation for right $\Dmod_X$-modules; it is defined at the level of 
derived categories by the formula
$$f_+ : {\bf D} (\Dmod_X) \longrightarrow {\bf D} (\Dmod_Y), \,\,\,\,\, \Nmod^{\bullet} \mapsto 
\derR f_* \big(\Nmod^{\bullet} \overset{\derL}{\otimes}_{\Dmod_X} \Dmod_{X\to Y} \big).$$
See \cite[\S1.5]{HTT} for more details, where this functor is denoted by $\int_f$.

Given a proper morphism of smooth varieties $f \colon X \rightarrow Y$, Saito has also constructed in \cite[\S2.3]{Saito-MHP} 
a filtered direct image functor 
$$f_+ : \DD^b \big({\rm FM}(\Dmod_X)\big) \rightarrow \DD^b \big({\rm FM}(\Dmod_Y)\big).$$
Here the two categories are the bounded derived categories of filtered $\Dmod$-modules on $X$ and $Y$ respectively.
Without filtration, it is precisely the functor above; some details about the filtration appear below, for the special 
$\Dmod$-modules considered here.

\subsection{Background on Hodge modules}\label{Hmod}
We now move to a review of the main notions and results from the theory of mixed Hodge modules that are used in this paper.

\noindent
{\bf Hodge modules and variations of Hodge structure.}
Let $X$ be a smooth complex variety of dimension $n$, and let $Z$ be an irreducible closed subset. 
Let $\V = (\mathcal{V}, F^\bullet, \V_{\QQ})$ be a polarizable variation of $\QQ$-Hodge structure of weight $\ell$ on an open set $U$ in the smooth locus of $Z$, where $(\shV, \nabla)$ is a vector bundle with flat connection with underlying $\QQ$-local system $\V_{\QQ}$, and Hodge filtration $F^\bullet \shV$.
Following \cite{Saito-MHP}, one can think of it as being smooth pure Hodge module of weight $\dim Z + \ell$ on U, whose main constituents are:

\begin{enumerate}
\item The left $\Dmod_U$-module $\Mmod = \mathcal{V}$ with filtration
$F_p \Mmod =  F^{-p} \mathcal{V}$.
\item The $\QQ$-perverse sheaf $P = \V_{\QQ} [n]$. 
\end{enumerate}

According to Saito's theory, this extends uniquely to a \emph{pure polarizable Hodge module} $M$ of weight $\dim Z + \ell$ on $X$, whose support is $Z$. This has an underlying perverse sheaf, which is the (shifted) intersection complex ${\rm IC}_Z( \V_{\QQ}) [n]=  {}^pj_{!*}\V_{\QQ}[n]$ associated to the given local system. 
It also has an underlying $\Dmod_X$-module, namely the minimal extension of $\Mmod$, corresponding to ${\rm IC}_Z( \V_{\CC})$
via the Riemann-Hilbert correspondence. Its filtration is (nontrivially) determined by the Hodge filtration on $U$; see e.g. Example \ref{normal_crossings} for the case of simple normal crossings boundary. 

More generally, in \cite{Saito-MHP} Saito introduced an abelian category ${\rm HM} (X, \ell)$ of pure polarizable Hodge modules on $X$ of weight $\ell$. The main two constituents of one such Hodge module $M$ are still:

\begin{enumerate}
\item A filtered (regular holonomic) left $\Dmod_X$-module $(\Mmod, F)$,  where $F = F_{\bullet} \Mmod$
is a good filtration by $\OX$-coherent subsheaves,  so that $\gr_{\bullet}^F \!
\Mmod$ is coherent over $\gr_{\bullet}^F \Dmod_X$. 

\item A $\QQ$-perverse sheaf $P$ on $X$ whose complexification corresponds to $\Mmod$ via the Riemann-Hilbert correspondence, so that there is an isomorphism
\[
	\alpha: \DR_X(\Mmod) \overset{\simeq}{\longrightarrow} P \tensor_{\QQ} \CC.
\]
\end{enumerate}

These are subject to a list of conditions, which are defined by induction on the dimension of the support of $M$. 
If $X$ is a point, a pure Hodge module is simply a polarizable Hodge structure of weight $\ell$. In general it is  required, roughly speaking, that the nearby and vanishing cycles associated to $M$ with respect to any locally defined holomorphic function are again Hodge modules, now on a variety of smaller dimension; for a nice discussion see \cite[\S12]{Schnell1}.
The definition of a polarization on $M$ is quite involved, but at the very least it involves an isomorphism ${\bf D} P \simeq P( \ell)$, 
where ${\bf D} P$ is the Verdier dual of the perverse sheaf $P$, together of course with further properties 
compatible with the inductive definition of Hodge modules.

Furthermore, M. Saito introduced in \cite{Saito-MHM} the abelian category $\MHM (X)$ 
of (graded-polarizable) \emph{mixed Hodge modules} on $X$. In addition to data as in (1) and (2) above, 
in this case a third main constituent is:

(3) A finite increasing weight filtration $W_{\bullet} M$ of $M$ by
objects of the same kind, compatible with $\alpha$, such that the graded quotients
$\gr_{\ell}^W M = W_{\ell} M / W_{\ell-1} M$ are pure Hodge modules in ${\rm HM} (X, \ell)$. 

Again, a mixed Hodge module on a point is a graded-polarizable mixed Hodge structure, while in general
these components are subject to several conditions defined by induction on the dimension of the support of $M$, involving
the graded quotients of the nearby and vanishing cycles of $M$; for a further discussion of the definition see 
\cite[\S20]{Schnell1}.  Only one class of examples of mixed (as opposed to pure) Hodge modules, described below, 
will appear explicitly in this paper. 

Let $D$ be a divisor in $X$ with complement $U$, and assume that we are given a variation of Hodge structure $\V$ on $U$.
Besides the pure Hodge module extension to $X$ described above,  it is also natural to consider a mixed Hodge module 
extension, denoted $j_*j^{-1} M$ in \cite{Saito-MHM}. Its underlying perverse sheaf is simply the direct image $j_*\V_\QQ$.
More precisely, 
\[ j_*j^{-1}M = \big( (\cV(*D), F); j_*\V_\QQ\big),\]
where $\cV(*D)$ is Deligne's meromorphic connection extending $\cV$ (see e.g. \cite[\S5.2]{HTT}). Further details are given in Example \ref{localization}.

\noindent
{\bf Strictness.}
From the point of view of applications, we mainly think here of Hodge modules as being filtered $\Dmod$-modules with extra 
structure. Part of this extra structure is the \emph{strictness} of the filtration; this property of $\Dmod$-modules underlying 
Hodge modules naturally appearing in birational geometry is crucial for applications.

Concretely, a morphism of filtered $\Dmod_X$-modules $f \colon (\Mmod, F) \rightarrow  (\Nmod, F)$, compatible with the filtrations, is called \emph{strict} if
$$f( F_k \Mmod)  = F_k \Nmod \cap f(\Mmod) \,\,\,\,\,\, {\rm for ~all~}k.$$ 
A complex $(\Mmod^\bullet, F_{\bullet} \Mmod^\bullet)$ is strict if all of its differentials are strict. 
This condition is equivalent to the injectivity, for every $i, k \in \ZZ$, of the induced morphism 
$$H^i (F_k \Mmod^\bullet) \longrightarrow H^i \Mmod^\bullet.$$
It is only in this case that the cohomologies of $\Mmod^\bullet$  can also be seen as filtered $\Dmod_X$-modules.
Moreover, if we denote by ${\rm FM}(\Dmod_X)$ the category of filtered $\Dmod_X$-modules, one can construct an associated bounded derived category ${\bf D}^b \big({\rm FM}(\Dmod_X)\big)$, and it continues to make sense to talk about the strictness of an object in ${\bf D}^b \big({\rm FM}(\Dmod_X)\big)$.

Recall now that  given a proper morphism of smooth varieties $f\colon X \rightarrow Y$,  we have a filtered direct image functor 
$$f_+ : \DD^b \big({\rm FM}(\Dmod_X)\big) \rightarrow \DD^b \big({\rm FM}(\Dmod_Y)\big).$$
For a filtered right $\Dmod_X$-module $(\Mmod, F)$, the strictness of  $f_+ (\Mmod, F)$ as an object in $\DD^b \big({\rm FM}(\Dmod_Y)\big)$ is equivalent  to the injectivity of the mapping
$$R^i f_* \big(F_k (\Mmod \overset{\derL}{\otimes}_{\Dmod_X} \Dmod_{X\to Y}) \big) \rightarrow 
R^i f_* (\Mmod \overset{\derL}{\otimes}_{\Dmod_X} \Dmod_{X\to Y})$$
for all integers $i$ and $k$, where the filtration on the left is the natural tensor product filtration induced by those on $\Mmod$ and 
$\Dmod_{X\to Y}$.

\noindent
{\bf Some fundamental theorems of Saito.}
It is time to state the most important results obtained by M. Saito in the theory of pure Hodge modules. They explain 
why Hodge modules are a natural and useful extension of the notion of variation of Hodge structure in the presence of singularities.

The first theorem describes the simple objects in the category of pure polarizable Hodge modules on $X$.  A 
pure Hodge module supported precisely along a subvariety $Z$ is said to have \emph{strict support} if it has no nontrivial subobjects or quotient objects whose support is $Z$.

\begin{theorem}[{\bf Simple objects}, {\cite[Theorem 3.21]{Saito-MHM}}]\label{structure}
Let $X$ be a smooth complex variety, and $Z$ an irreducible closed subvariety of $X$. Then:
\begin{enumerate}
\item Every polarizable VHS of weight $\ell$ defined on a nonempty open set of $Z$ extends uniquely 
to a pure polarizable Hodge module of weight $\ell + \dim Z$ with strict support $Z$.
\item Conversely, every pure polarizable Hodge module of weight $\ell + \dim Z$ with strict support $Z$ is obtained in this way.
\end{enumerate}
\end{theorem}

The category ${\rm HM} (X, \ell)$ is in fact semi-simple, due to the existence of polarizations: each object admits a 
decomposition by support, and simple objects with support equal to an irreducible subvariety $Z$ are precisely those in the theorem.

One of the most important results of the theory is Saito's theorem on the behavior 
of direct images of polarizable Hodge modules via projective morphisms, only a part of which is stated below.

\begin{theorem}[{\bf Stability Theorem}, {\cite[Th\'eor\`eme 5.3.1]{Saito-MHP}}]\label{stability}
Let $f\colon X \rightarrow Y$ be a projective morphism of smooth complex varieties, and $M \in {\rm HM} (X, \ell)$  a polarizable pure Hodge module on $X$ with underlying filtered \emph{right} $\Dmod$-module $(\Mmod, F)$.
Then the filtered direct image $f_+ (\Mmod, F)$ is strict, and $\mathcal{H}^i f_+ \Mmod$ underlies a polarizable 
Hodge module $M_i  \in {\rm HM}(Y, \ell + i)$.
\end{theorem}

The strictness of the direct image is a key property of $\Dmod$-modules underlying Hodge modules that is not shared by arbitrary filtered 
$\Dmod$-modules. It is one of the most crucial inputs for the applications described in this survey, by means of the following 
consequences:

\noindent
(i) (\emph{Laumon's formula}, see e.g. \cite[\S2.2, Theorem 9]{mhmgv} for a proof.) If $\shA_X : = {\rm Sym}~T_X \simeq \gr^F_{\bullet} \Dmod_X$, and $\shA_Y$ is defined similarly on $Y$, then
$$\gr^F_{\bullet} \mathcal{H}^i f_+ (\Mmod, F) \simeq R^i f_* \big(\gr^F_{\bullet} \Mmod \otimes_{\shA_X} f^* \shA_Y \big),$$ 

\noindent
(ii) (\emph{Saito's formula}, see \cite[2.3.7]{Saito-MHP}.) The associated graded of the filtered de Rham complex satisfies
$$R^i f_*  \gr_k^F {\rm DR} (\Mmod, F) \simeq \gr_k^F {\rm DR} \big(\mathcal{H}^i f_+ (\Mmod, F)\big) .$$
 
\noindent
(iii) (\emph{Degeneration of the Hodge-to-de Rham spectral sequence}.)
In the case of the constant morphism $f \colon X \to {\rm pt}$, the strictness of $f_+ (\Mmod, F)$ is equivalent to the $E_1$-degeneration 
of the  Hodge-to-de Rham spectral sequence  
$$E_1^{p,q} = \mathbf{H}^{p+ q}  \big(X, \gr_{-q}^F \DR ( \Mmod) \big) \implies  H^{p+q} \big(X, \DR ( \Mmod)\big).$$

Finally, a fundamental consequence of (the full version of) the Stability Theorem is the extension of the well-known Decomposition Theorem in topology to the setting of pure polarizable Hodge modules \cite{Saito-MHP}; here I state the filtered $\Dmod$-modules version, used later.

\begin{theorem}[{\bf Saito Decomposition Theorem}]\label{decomposition}
Let $f\colon X \to Y$ be a projective morphism of smooth complex varieties, 
and let $M \in {\rm HM} (X, \ell)$, with underlying filtered $\Dmod$-module $(\Mmod, F)$. Then 
$$ f_+(\Mmod, F) \simeq \bigoplus_{i\in \ZZ} \mathcal{H}^i f_+ (\Mmod, F) [-i]$$
in the derived category of filtered $\Dmod_Y$-modules.
\end{theorem}

\noindent
{\bf Examples.}
The Hodge modules that have proved most useful in birational geometry to date are among the simplest concrete examples of such objects
that one can write down; nevertheless, they contain deep information.

 \begin{example}[{\bf The trivial Hodge module}]\label{trivial_HM}
If $X$ is smooth of dimension $n$ and $\V = \QQ_X$ is the constant variation of Hodge structure, we have that $P = \QQ_X[n]$,  while
$\Mmod = \shO_X$ with the natural left $\Dmod_X$-module structure and trivial filtration, so that $\gr_k^F \shO_X$ is $\shO_X$ for 
$k = 0$ and $0$ for $k \neq 0$. The associated Hodge module is usually denoted $\QQ_X^H[n]$. The de Rham complex of $\Mmod$ is 
\[
 \DR_X (\shO_X)  = \Bigl\lbrack
		\shO_X \to \Omega_X^1  \to \dotsb \to \Omega_X^n 
	\Bigr\rbrack [n],
\]
i.e. the usual de Rham complex of $X$ considered in degrees $-n, \ldots, 0$. With the filtration induced from that on $\shO_X$, 
note that 
$$\gr_{-k}^F \DR (\shO_X) = \Omega^k_X [n-k]\,\,\,\,  
{\rm ~for~ all~}k.$$
In this example we have $p(\Mmod) = 0$ and $S (\Mmod) = \shO_X$.
\end{example}

\begin{example}[{\bf Direct images}]\label{direct_image}
Let $f \colon X \rightarrow Y$ be a projective morphism of smooth complex varieties of dimensions $n$ and $m$ respectively, 
and let $\V$ be a polarizable variation of $\QQ$-Hodge structure of weight $k$ on an open dense subset 
$U \subset X$, 
inducing a pure Hodge module $M$ of weight  $n + k$ with strict support $X$. Here it is convenient to use right
$\Dmod$-modules: if $(\Nmod, F)$ is  filtered right $\Dmod_X$-module underlying $M$, Theorem \ref{decomposition} gives a decomposition
$$f_+ (\Nmod, F) \simeq \bigoplus_i (\Nmod_i, F) [-i]$$ 
in the derived category of filtered  $\Dmod_Y$-modules.
According to Theorem \ref{stability}, each $(\Nmod_i, F)$ underlies a 
 pure Hodge module $M_i = \mathcal{H}^i f_*M$ on $Y$, of weight $n+k +i$. 
Furthermore, $f_+ (\Nmod, F)$ is strict; using Saito's formula (ii) in the previous subsection this implies for instance that
\begin{equation}\label{lowest_decomposition}
\derR f_* S (\Nmod) \simeq F_{p(\Nmod)} (f_+ \Nmod) \simeq \bigoplus_i F_{p(\Nmod)} \Nmod_i [ -i]
\end{equation}
in the bounded derived category of coherent sheaves on $Y$. 

For example,  when $\V$ is the constant variation of Hodge structure $\QQ_X$, by Example \ref{trivial_HM} we have
$p (\Nmod) = -n$ and $S(\Nmod) = \omega_X$. This implies for all $i$ that
$$p(N_i) = m-n \,\,\,\, {\rm and} \,\,\,\, F_{m-n} \Nmod_i = R^i f_* \omega_X.$$
Consequently, formula (\ref{lowest_decomposition}) specializes to 
$$\derR f_* \omega_X \simeq \bigoplus_i R^i f_* \omega_X [ -i],$$
which is the well-known Koll\'ar decomposition theorem \cite{Kollar3}.\footnote{To recover the full statement of Koll\'ar's theorem 
note that we should allow $Y$ to be singular, and we can indeed do this. One can work with Hodge modules on singular varieties 
using local embeddings into smooth varieties and Kashiwara's theorem; for details see \cite[\S2.1]{Saito-MHM} or 
\cite[\S7]{Schnell2}.}
In fact $R^i f_* S(\Nmod)$ satisfy other important properties known from \cite{Kollar1} in the case of canonical bundles; for instance, 
answering a conjecture of Koll\'ar, Saito proved the following:
 
 \begin{theorem}[Saito, \cite{Saito-Kollar}]\label{saito_GR}
 Let $f \colon X \rightarrow Y$ be  a surjective projective morphism, with $X$ smooth, and let $(\Mmod, F)$ be the 
filtered $\Dmod$-module underlying a pure Hodge module with strict support $X$ that is generically a polarized variation of Hodge structure $\V$. For each $i\ge 0$, one has 
 $$R^i f_* S(\Mmod) = S(Y, \V^i),$$
 the lowest Hodge piece of the variation of Hodge structure $\V^i$ on the intersection cohomology of $\V$ along
 the fibers of $f$. Consequently, $R^i f_* S(\Mmod)$ are torsion-free, and in particular 
 $$R^i f_* S(\Mmod) = 0 \,\,\,\,\,\,{\rm for} \,\,\,\, i > \dim X - \dim Y.$$
 \end{theorem}
\end{example}

\begin{example}[{\bf Localization}]\label{localization}
Let $\Mmod$ be a left $\Dmod_X$-module and $D$ an effective divisor on a smooth variety $X$, given locally by an equation $f$. 
One can define a new $\Dmod_X$-module $\Mmod (*D)$ by localizing $\Mmod$ in $f$; in other words, globally we have
$$\Mmod(*D) =  j_* j^{-1} \Mmod,$$
where $j: U \hookrightarrow X$ is the inclusion of the complement $U = X \smallsetminus D$.

When $\Mmod$ underlies a mixed Hodge module $M$,  $\Mmod (*D)$ underlies the corresponding mixed Hodge module $j_* j^{-1} M$, and so continues to carry a natural Hodge filtration $F$. This is in general very complicated to compute; for instance, 
the case $\Mmod = \shO_X$, where $\shO_X (*D)$ is the sheaf of meromorphic functions on $X$ that are holomorphic on $U$ and the corresponding Hodge module is $j_* \QQ^H_U \decal{n}$, is very subtly related to the singularities of $D$, and to Deligne's Hodge filtration 
on $H^{\bullet} (U, \CC)$. It is known for example that
$$S \big( \shO_X (*D) \big) = \shO_X (D) \otimes \mathcal{I} \big(X, (1-\epsilon) D\big),$$
where ideal in question is the multiplier ideal of the $\QQ$-divisor $(1-\epsilon)D$, with $0 < \epsilon \ll 1$.
Much more on the Hodge filtration on $\shO_X (*D)$ is discussed in Section \ref{Hodge_ideals}.
\end{example}

\begin{example}[{\bf Normal crossings case}]\label{normal_crossings}
Often one analyzes the Hodge filtration on a Hodge $\Dmod$-module of interest by pushing forward a better understood filtration 
on a log resolution, and applying the Stability Theorem. This example shows that for both the pure and mixed Hodge module 
extensions of a VHS over simple normal crossings boundary, there are explicit formulas for the underlying filtered $\Dmod$-module.

Let $(X, D)$ be a pair consisting of a smooth projective variety $X$ of  dimension $n$, and a simple normal crossings divisor $D$.
Denote $U=X \smallsetminus  D$ and $j: U\hookrightarrow X$. We consider a polarizable VHS
\[\V=(\cV, F_{\bullet}, \V_\QQ)\]
over $U$, with quasi-unipotent local monodromies along the components $D_i$ of $D$. In particular the eigenvalues of all residues are rational numbers. For $\alpha \in \ZZ$,  we 
denote by $\cV^{\geq\alpha}$ (resp. $\cV^{>\alpha}$) the Deligne extension with eigenvalues of residues along the  $D_i$ in $[\alpha, \alpha+1)$ (resp. $(\alpha, \alpha+1]$). Recall that $\cV^{\geq\alpha}$ is filtered by
\begin{equation}\label{filtr}
F_p\cV^{\geq\alpha}=\cV^{\geq\alpha}\cap j_*F_p\cV,
\end{equation}
while the filtration on $\cV^{>\alpha}$ is defined similarly. The terms in the filtration are 
locally free by Schmid's nilpotent orbit theorem \cite{Schmid} (see also e.g. \cite[2.5(iii)]{Kollar3} for the quasi-unipotent case).

If $M$ is the pure Hodge module with strict support $X$ uniquely extending 
$\V$, it is proved in \cite[\S 3.b]{Saito-MHM}  that
\[M=\big(\Dmod_X\cV^{>-1}, F_{\bullet}, j_{!*}\V_\QQ\big),\] where
\[F_p\Dmod_X\cV^{>-1}=\sum_i F_i\Dmod_X \cdot F_{p-i}\cV^{>-1}.\]
For the natural mixed Hodge module extension of $\V$, namely 
\[j_*j^{-1}M =\big(\cV(*D), F_{\bullet}, j_*\V_\QQ\big),\]
recall that $\cV(*D)$ is Deligne's meromorphic connection extending $\cV$. It has a lattice 
defined by $\cV^{\geq\alpha}$ for any $\alpha\in\QQ$, i.e. $\cV(*D)=\cV^{\ge \alpha}\otimes \shO_X (*D)$,
and its filtration is given by
\[F_p\cV(*D)=\sum_i F_i\Dmod_X \cdot F_{p-i}\cV^{\geq-1}.\]
\end{example}

\section{Vanishing and positivity package for Hodge modules}\label{vanishing_positivity}

One of the key aspects regarding the applications of Hodge module theory to birational geometry is the existence 
of a full package of vanishing and positivity theorems generalizing well-known results about (direct images of) canonical bundles; many of those had of course themselves been obtained using more classical methods in Hodge theory.

\subsection{Vanishing theorems}\label{vanishing}
We start with a series of vanishing results for mixed Hodge modules, which includes generalizations of most of the 
standard vanishing theorems.

\noindent
{\bf Kodaira-Saito vanishing.}
Saito noted in his original paper \cite{Saito-MHM} that the Kodaira-Nakano vanishing theorem has a far-reaching extension to the setting of mixed Hodge modules. 

\begin{theorem}[{\cite[\S2.g]{Saito-MHM}}]\label{saito_vanishing}
Let $X$ be a complex projective variety, and $L$ an ample line bundle on $X$. Consider an integer $m >0 $ such that 
$L^{\otimes m}$ is very ample 
and gives an embedding $X \subseteq \PP^N$. Let $(\Mmod, F)$ be the filtered $\Dmod$-module underlying a mixed Hodge module $M$ on 
$\PP^N$ with support contained in $X$. Then:

\noindent
(1) $\gr_k^F \DR_{\PP^N} (\Mmod)$ is an object in $\DD^b (X)$ for each $k$, independent of the embedding of $X$ in $\PP^N$.

\noindent
(2) We have the hypercohomology vanishing 
$${\bf H}^i \bigl( X, \gr_k^F \DR_{\PP^N} (\Mmod) \tensor L \bigr) = 0 {\rm ~ for~ all~} i > 0.$$
and 
$${\bf H}^i \bigl( X, \gr_k^F \DR_{\PP^N} (\Mmod) \tensor L^{-1} \bigr) = 0\,\,\,\,{\rm~ for~ all~} i < 0.$$
\end{theorem}

\begin{remark}
The sightly awkward formulation in the statement is due to the fact that $X$ is allowed to be singular, which is important in applications 
(see for instance Example \ref{classical_vanishing}(ii) below). See \cite[\S7]{Schnell2} and \cite[\S8]{Popa} for more on this. 
When $X$ is smooth, one can of course work directly with a mixed Hodge module on $X$. 
\end{remark}

The proof given in \cite{Saito-MHM} goes along the lines of Ramanujam's topological approach to Kodaira vanishing; a detailed 
account is given in \cite{Popa}. Another proof, this time along the lines of an approach to vanishing theorems due to 
Esnault-Viehweg, was given in \cite{Schnell2}.

\begin{example}\label{classical_vanishing}
(i) ({\bf Kodaira-Nakano vanishing.}) Let $X$ be a smooth projective complex variety of dimension $n$.
For the trivial Hodge module $M = \QQ_X^H \decal{n}$, according to Example \ref{trivial_HM} 
the corresponding left $\Dmod$-module is $\shO_X$, and 
$$\gr_{-k}^F \DR_X(\shO_X) = \Omega_X^k \decal{n-k} \,\,\,\, {\rm for ~all~} k.$$ 
In this case Theorem \ref{saito_vanishing} is therefore the Kodaira-Nakano vanishing theorem.

As this is used repeatedly, it is worth recording the fact that the Kodaira vanishing theorem, which corresponds to the 
lowest non-zero graded piece of the filtration on $\omega_X$, can be seen as an example of the simplest special case of Theorem \ref{saito_vanishing}.

\begin{corollary}\label{lowest}
If $(\Mmod, F)$ is the filtered $\Dmod$-module underlying a mixed Hodge module $M$ on a 
projective variety $X$,\footnote{In the sense of Theorem \ref{saito_vanishing} if $X$ is singular.} and $L$ is an ample line bundle on $X$, then
$$H^i \big(X,  \omega_X \otimes S(\Mmod) \otimes L \big) = 0 \,\,\,\,{\rm for~all~} i > 0.$$
\end{corollary}

\noindent
(ii)  ({\bf Koll\'ar vanishing.}) Let $f \colon X \rightarrow Y$ be a morphism of complex projective varieties, with $X$  smooth of dimension $n$, and let $L$ be an ample line bundle on $Y$. Considering the push-forward $f_* \QQ_X^H \decal{n}$, according to Example \ref{direct_image} for the underlying \emph{right} $\Dmod$-module we have  
$$f_+ (\omega_X, F)  \simeq \bigoplus_i (\Nmod_i, F) [-i]$$ 
in the derived category of filtered $\Dmod_Y$-modules, and for 
each $i$ we have $S(\Nmod_i) = R^i f_* \omega_X$. Therefore here a special case of Corollary \ref{lowest} is Koll\'ar's vanishing theorem \cite[Theorem 2.1(iii)]{Kollar1}, saying that for any ample line bundle $L$ on $Y$, 
$$H^i (Y, R^i f_* \omega_X \otimes L ) = 0 \,\,\,\,\,\, {\rm ~for~all} \,\, i >0 \,\,{\rm and ~all~} j.$$

\noindent 
(iii) ({\bf Nadel vanishing.}) If $D$ is an effective divisor on a smooth projective variety $X$, recall from Example \ref{localization} that 
$\shO_X(*D)$ underlies a mixed Hodge module, and that its lowest term in the Hodge filtration is
$\shO_X (D) \otimes \mathcal{I} \big(X, (1-\epsilon) D\big)$. Using Corollary \ref{lowest}, it follows that 
$$H^i \big(X, \omega_X\otimes L \otimes \mathcal{I} \big(X, (1-\epsilon) D\big) \big) = 0 \,\,\,\,\,\,{\rm for ~all} \,\,\,\,i > 0,$$
for any line bundle $L$ such that $L (-D)$ is ample. This is a special case of the Nadel vanishing theorem, see \cite[Theorem 9.4.8]{Lazarsfeld}. Nadel vanishing for arbitrary $\QQ$-divisors cannot be seen directly as an example of Saito vanishing, but can be deduced from it 
using covering constructions.
\end{example}

\noindent
{\bf Kawamata-Viehweg-type vanishing.}
In \cite[\S11]{Popa} it was noted that a vanishing theorem for big and nef divisors, analogous to the well-known Kawamata-Viehweg theorem, holds for the lowest term $S(\Mmod)$ of the Hodge filtration, under a transversality hypothesis involving the non-smooth locus of $M$ and the augmented base locus of $L$. This hypothesis has been removed by Suh and Wu independently.

\begin{theorem}[{\cite{Suh}, \cite{Wu}}]\label{KV}
Let $X$ be a complex projective variety, $L$ a big and nef line bundle on $X$, 
and $(\Mmod, F)$ the filtered left $\Dmod$-module underlying a polarizable pure Hodge module with strict support $X$.
Then
$$H^i \big(X, \omega_X \otimes S(\Mmod) \otimes L\big) = 0 \,\,\,\,\,\, {\rm for ~all} \,\,\,\, i > 0.$$
\end{theorem}

The approach in \cite{Suh} is in the spirit of the Esnault-Viehweg method using the degeneration of the Hodge-to-de Rham spectral sequence. In the process Suh obtains a related statement for Deligne canonical extensions that applies to all graded quotients of the Hodge filtration; this itself has interesting applications, see Ch.\ref{other_applications}. The approach in \cite{Wu}  uses an induction based on adjunction-type formulas involving mixed Hodge modules of nearby cycles, more in line with Kawamata's approach to vanishing.

\begin{remark}
One cannot expect the statement in Theorem \ref{KV} to hold for other graded pieces $\gr_k^F \DR (\Mmod)$ with $k \neq p(\Mmod)$, 
as in Theorem \ref{saito_vanishing}. Indeed, in the case $M = \QQ_X^H [n]$ we have seen in Example \ref{classical_vanishing} that these graded pieces are (shifts of) various $\Omega_X^p$; it is well-known however that the analogue of Nakano vanishing does not hold for big and nef divisors. 
\end{remark}

\noindent
{\bf Injectivity.}
On a related note, Wu \cite{Wu} has shown that the vanishing results involving the lowest term $S(\Mmod)$ of the Hodge filtration admit a further improvement extending Koll\'ar's injectivity theorem \cite[Theorem 2.2]{Kollar1} and its generalizations by 
Esnault-Viehweg \cite[5.1]{EV}. The proof relies on Suh's approach to vanishing for Deligne canonical extensions in \cite{Suh}.

\begin{theorem}[{\cite[Theorem 1.4]{Wu}}]\label{injectivity}
Let $X$ be a complex projective variety, $D$ an effective divisor on $X$, 
and $(\Mmod, F)$ the filtered left $\Dmod$-module underlying a polarizable pure Hodge module with strict support $X$.
Let $L$ be a line bundle on $X$ satisfying either one of the following properties:

\begin{enumerate}
\item
$L$ is nef and big

\item
$L$ is semiample, and $H^0 \big(X, L^{\otimes m} (- D) \big) \neq 0$ for some $m > 0$.
\end{enumerate}
Then the natural mapping 
$$H^i \big(X, \omega_X \otimes S(\Mmod) \otimes L \big) \longrightarrow H^i \big(X, \omega_X \otimes S(\Mmod) \otimes L (D) \big) $$
is injective for all $i$. 
\end{theorem}

\noindent
{\bf $\QQ$-divisors.}
Just as in the usual setting, in view of future applications it is important to have extensions of the results above to the 
case of $\QQ$-divisors. In \cite[Theorem 11.1]{Popa}, a Kawamata-Viehweg-type vanishing theorem for arbitrary $\QQ$-divisors is proved under the transversality hypothesis  mentioned above.
On the other hand, \cite[\S7]{Wu} contains $\QQ$-versions of Theorem \ref{KV} and Theorem \ref{injectivity} with no such hypotheses;  these necessarily require extra conditions on the $\QQ$-coefficients in the abstract Hodge theoretic setting, but do recover the standard vanishing theorems for $\QQ$-divisors. Recently Wu \cite{Wu2} has announced a further generalization of the Kawamata-Viehweg vanishing theorem to the setting of Hodge modules, involving arbitrary $\QQ$-divisors.

\subsection{Weak positivity}\label{scn:WP}
A general property of the lowest term of the Hodge filtration of a generically defined polarizable variation of (mixed) Hodge structure, or of an extension of such over simple normal crossings boundary, is its semipositivity. This was studied in a celebrated series of works by Griffiths, Fujita, Kawamata, Koll\'ar and others; for the most recent results in this direction, and detailed references, see Fujino-Fujisawa-Saito \cite{FFS}.

To deal with more general situations, Viehweg introduced a more relaxed notion of positivity, still sufficient for many applications: a torsion-free coherent sheaf $\shF$ on a quasiprojective variety $X$ is \emph{weakly positive} on a non-empty open set $U \subseteq X$ if for every ample line bundle $A$ on $X$ and every $a\in \NN$, the sheaf ${\hat S}^{ab} \shF \otimes A^{\otimes b}$ is generated by global sections at each point of $U$ for $b$ sufficiently large.\footnote{Here ${\hat S}^p \shF$ denotes the reflexive hull of the symmetric power $S^p \shF$.}
For line bundles for instance, while semipositivity is analogous to nefness in algebraic language, weak positivity is equivalent to pseudo-effectivity.
Viehweg \cite{Viehweg1} showed for instance that for any morphism of smooth projective varieties $f\colon Y \to X$, and any $m \ge 1$, the 
sheaf $f_* \omega_{Y/X}^{\otimes m}$ is weakly positive.

These results have useful versions at the level of Hodge modules. For the rest of the section, we fix 
a polarizable pure Hodge module $M$ with strict support $X$, and consider the filtered left $\Dmod_X$-module 
$(\Mmod, F)$ underlying $M$. It is not hard to see that $S(\Mmod)$ is a torsion-free sheaf. 
First, along the lines of results regarding the lowest term in the Hodge filtration, one has the following extension to this setting:

\begin{theorem}[{\cite[Theorem 1.4]{Schnell3}, \cite[Theorem 10.4]{Popa}}]\label{WP_old}
The sheaf $S(\Mmod)$ is weakly positive torsion-free sheaf.
\end{theorem} 

\begin{example}
When $f\colon Y \to X$ is a morphism of smooth projective varieties, taking $M = M_0$ in Example \ref{lowest_decomposition}, 
one recovers Viehweg's result above for $m = 1$, i.e. the fact that $f_* \omega_{Y/X}$ is weakly positive. Moreover, taking $M = M_i$, 
one gets the same result for $R^i f_* \omega_{Y/X}$ for all $i$.

A sketch of the proof of Theorem \ref{WP_old} may be instructive 
here, as it shows the efficiency of considering Hodge modules, as well as an approach originally due to Koll\'ar using vanishing; it provides a quicker argument even in the setting of Viehweg's theorem. In brief, considering for any $m>0$ the diagonal embedding 
$$i \colon X \hookrightarrow X\times \cdots \times X,$$
where the product is taken $m$ times, and pulling back the box-product $M\boxtimes \cdots \boxtimes M$ via $i$, 
it is not too hard to show that there exists an abstract Hodge module $M'$ on $X$, together with a morphism
$$S(\Mmod^\prime) \longrightarrow S(\Mmod)^{\otimes m},$$
which is an isomorphism on the open set where $M$ is a VHS. Now taking any very ample line bundle $L$ on $X$, 
Corollary \ref{lowest} combined with Castelnuovo-Mumford regularity shows that 
$S(\Mmod^\prime) \otimes \omega_X \otimes L^{\otimes (n+1)}$ is globally generated, where $n = \dim X$. Since $m$ is arbitrary, 
a standard argument implies the weak positivity of $S(\Mmod)$.
\end{example}

It turns out that there is a more general result, which holds at each step of the Hodge filtration; this is very useful 
for applications. It has its origin in previous results of Zuo \cite{Zuo} and Brunebarbe \cite{Brunebarbe}, which are 
analogous Griffiths-type metric statements in the setting of Deligne canonical extensions over a simple normal crossings boundary.
Concretely, for each $p$ we have a natural Kodaira-Spencer type $\shO_X$-module homomorphism
$$\theta_p:  {\rm gr}_p^F \Mmod \longrightarrow {\rm gr}_{p+1}^F \Mmod \otimes \Omega_X^1$$
induced by the $\Dmod_X$-module structure, and we denote 
$$K_p (\Mmod) : = {\rm ker}~ \theta_p.$$

\begin{theorem}[{\cite[Theorem A]{PW}}]\label{thm:WP}
The torsion-free sheaf $K_p (\Mmod)^\vee$ is weakly positive for any $p$.
\end{theorem}

One can show that Theorem \ref{thm:WP} implies Theorem \ref{WP_old} via duality arguments. On the other hand, the proof of Theorem \ref{thm:WP} 
is done by reduction to the metric statements in the papers cited above; it would be interesting to know whether it has a proof based on vanishing theorems as well.

\section{Generic vanishing}\label{generic_vanishing}

Hodge modules have made an impact in generic vanishing theory. Besides theoretical results on the graded 
pieces of the Hodge filtration that are useful in applications, as we will see later on, their study has led to
previously inaccessible results on bundles of holomorphic forms (in the projective setting) and higher direct images of canonical bundles, as well as structural results (in the K\"ahler setting).

Let $T$ be a compact complex torus. Given a coherent sheaf $\shF$ on $T$, and an integer $i\ge 1$, one defines 
\[
	V^i (\shF) = \menge{L \in \Pic^0(T)}{H^i(T, \shF \otimes L) \neq 0}.
\]
Recall the following definitions from \cite{PP3} and \cite{PP4} respectively; the sheaf $\shF$ is called a \define{GV-sheaf} if the inequality
\[
	\codim_{\Pic^0(T)} V^i(\shF) \geq i
\]
is satisfied for every integer $i \geq 0$. It is called \define{M-regular} if the inequality
\[
	\codim_{\Pic^0(T)} V^i(\shF) \geq i+1
\]
is satisfied for every integer $i \geq 1$. The notion of $GV$-sheaf is meant to formalize the behavior of canonical bundles of 
compact K\"ahler manifolds with generically finite Albanese maps, which were shown to have this property in \cite{GL1}. 
Using vanishing theorems, Hacon \cite[\S4]{Hacon} showed the stronger fact that if $f \colon X \to A$ is a morphism from a smooth projective variety to an abelian variety, then the higher direct image sheaves $R^j \fl \omX$ are $GV$-sheaves on $A$; 
see also \cite{PP3} for further developments. On the other hand, the stronger notion of $M$-regularity turns out to be very useful in the study of syzygies and of linear series on irregular varieties; see e.g \cite{PP5} and the references therein.

\subsection{Generic vanishing on smooth projective varieties}\label{scn:GV-projective}
Given a morphism $f \colon X \to A$ to an abelian variety, unlike in the case of $R^j f_* \omega_X$, with 
the techniques developed in the papers above it has not been possible to obtain optimal results regarding the $GV$-property of other bundles of holomorphic forms or, in an appropriate sense, of rank $1$ local systems. This was overcome in \cite{mhmgv} by a series of generic vanishing results involving mixed Hodge modules; here is a first version:

\begin{theorem}[{\cite[Theorem 1.1]{mhmgv}}]\label{hm_abelian}
Let $A$ be a complex abelian variety, and $M$ a mixed Hodge module on $A$ with
underlying filtered $\Dmod$-module $(\Mmod, F)$. Then, for each $k
\in \ZZ$, the coherent sheaf $\gr_k^F \Mmod$ is a $GV$-sheaf on $A$.
\end{theorem}

Note that by Examples \ref{trivial_HM} and \ref{direct_image}, this result includes the theorems of Green-Lazarsfeld and 
Hacon described above. For a subsequent generalization, see Theorem \ref{thm:Chen-Jiang} below. The main reason 
for Theorem \ref{hm_abelian}  is the fact that on abelian varieties Kodaira-Saito vanishing holds in a much stronger 
form than Theorem \ref{saito_vanishing}:

\begin{proposition}[{\cite[Lemma 2.5]{mhmgv}}]\label{abelian}
Under the assumptions of Theorem \ref{hm_abelian}, if $L$ is an ample line bundle, then 
\[
        H^i \bigl( A, \gr_k^F \Mmod \tensor L \bigr) = 0 \,\,\,\, {\rm for~all ~} i >0.
\]
\end{proposition}

This has a quite simple explanation, due to the triviality of the cotangent bundle: if $g = \dim A$, then we have 
\[
	\gr_k^F \DR_A(\Mmod) = \Big\lbrack
		\gr_k^F \Mmod \to \Omega_A^{1} \tensor \gr_{k+1}^F \Mmod \to \dotsb \to
			\Omega_A^{g} \tensor \gr_{k+g}^F \Mmod
	\Big\rbrack,
\]
supported in degrees $-g, \dotsc, 0$. According to Theorem~\ref{saito_vanishing},
this complex satisfies 
\[
	{\bf H}^i \bigl( A, \gr_k^F \DR_A(\Mmod) \otimes L \bigr) = 0 \,\,\,\, {\rm for~all ~} i >0.
\]
On the other hand, using the fact that $\Omega_A^{1} \simeq \shO_A^{\oplus g}$, and $\gr_k^F \Mmod = 0$ for $k \ll 0$,
one has for each $\ell$ a distinguished triangle 
\[	
	E_\ell \to \gr_\ell^F \DR_A(\Mmod) \to \gr_{\ell+g}^F \Mmod \to E_\ell [1],
\]
where it can be assumed inductively that $E_\ell$ is an object satisfying $H^i ( A, E_\ell \otimes L) = 0$. 
One can then deduce the statement of Proposition \ref{abelian} by induction on $k$. 
Theorem \ref{hm_abelian} can be  deduced from Proposition \ref{abelian} by using a cohomological criterion of Hacon \cite{Hacon} for detecting the $GV$-property. 

Theorem \ref{hm_abelian} implies an optimal Nakano-type generic vanishing theorem, by taking $(\Mmod, F)$ 
to be the various constituents of the push-forward $f_+ (\omega_X, F)$ via the Decomposition Theorem, as in Example \ref{direct_image}. Here $f\colon X \to A$ is the Albanese map of $X$. Recall that 
\define{defect of semismallness} of $f$ is defined by the formula
\[
	\delta(f) = \max_{\ell \in \NN} \bigl( 2 \ell - \dim X + \dim A_{\ell} \bigr),
\]
where $A_{\ell} = \menge{y \in A}{\dim f^{-1}(y) \geq \ell}$ for $\ell \in \NN$. The condition $\delta (f) = 0$ is equivalent to $f$ 
being semismall.

\begin{theorem}[{\cite[Theorem 1.2]{mhmgv}}]\label{thm:nakano}
Let $X$ be a smooth complex projective variety of dimension $n$. Then\footnote{Here, by analogy with the definition for sheaves
on tori, one sets $$V^q(\Omega_X^p) = \menge{L \in \Pic^0(X)}{H^q(X, \Omega_X^p \otimes L) \neq 0}.$$}
\[
 	\codim V^q(\Omega_X^p) \ge \abs{p + q - n} - \delta(f) \,\,\,\,\,\,{\rm for~all} \,\,\,\,p,q \in \NN, 
\]	
and there exist $p,q$ for which equality is attained. In particular, 
$f$ is semismall if and only if $X$ satisfies the generic Nakano vanishing theorem, i.e.
\[
	\codim V^q (\Omega_X^p)  \geq \abs{p + q - n}  \,\,\,\,\,\,{\rm for~all} \,\,\,\,p,q \in \NN.
\]
\end{theorem}

Working with Hodge modules allows one however to obtain a much more complex picture than Theorem \ref{hm_abelian},
putting in some sense topologically trivial line bundles and holomorphic $1$-forms on equal footing. 
If $(\Mmod, F)$ is the filtered $\Dmod$-module underlying a mixed Hodge module on $A$, 
we denote by $\gr^F \! \Mmod$ the coherent sheaf on $T^* A = A \times V$ corresponding to the total 
associated graded $\gr_\bullet^F \Mmod$, where $V = H^0 (A, \Omega_A^1)$.
If $P$ is the normalized Poincar\'e bundle on $A \times \Ah$, using its pull-back to
$A \times \Ah \times V$ as a kernel, we define the (relative) Fourier-Mukai transform of $\gr^F
\! \Mmod$ to be the complex of coherent sheaves on $\Ah \times V$ given by
\[
	E = \derR \Phi_P(\gr^F \! \Mmod) 
		= \derR (p_{23})_{\ast} \Bigl( p_{13}^{\ast} \bigl( \gr^F \! \Mmod \bigr) 
			\tensor p_{12}^{\ast} P \Bigr).
\]

\begin{theorem}[{\cite[Theorem 1.6]{mhmgv}}]\label{thm:mCoh}
With the notation above, let $(\Mmod, F)$ be the filtered left $\Dmod$-module associated to 
any of the direct summands of $f_+ (\omega_X, F)$, where 
$f \colon X \to A$ is a morphism from a smooth projective variety. Then the complex $E$ has the following 
properties:
\begin{enumerate}
\item 
$E$ is a perverse coherent sheaf, meaning that its cohomology sheaves
$\shH^{\ell} E$ are zero for $\ell < 0$, and 
$${\rm codim}~{\rm Supp}(\shH^{\ell} E) \ge 2\ell \,\,\,\,\,\,{\rm for} \,\,\,\,\ell \ge 0.$$
\item
The union of the supports of all the higher cohomology sheaves of $E$ is a finite
union of translates of triple tori in $\Ah \times V$.
\item 
The dual complex $\derR \shHom(E, \shO)$ has the same properties.
\end{enumerate}
\end{theorem}
Here a triple torus, in Simpson's terminology \cite{Simpson}, is a subset of the form 
\[
	\im \Bigl( \varphi^{\ast} \colon \hat{B} \times H^0(B, \OmB^1)
		\to \Ah \times H^0(A, \OmA^1) \Bigr),
\]
where $\varphi \colon A \to B$ is a morphism to another abelian variety. The properties of Hodge $\Dmod_A$-modules of 
geometric origin provided by the theorem are very useful in applications; for instance they are key to the proof of 
Theorem \ref{one-forms} below.

Theorem \ref{thm:mCoh} turns out to hold for arbitrary holonomic $\Dmod$-modules on $A$; this, as well as numerous 
other similar results for $\Dmod$-modules on abelian varieties, were obtained by Schnell  in \cite{Schnell}.  Moreover, similar results are obtained in \cite{mhmgv} and \cite{Schnell} involving a different type of Fourier-Mukai transform, namely 
the Laumon-Rothstein transform 
\[
	\derR \Phi_{P^\natural} : \DD_{{\rm coh}}^b (\Dmod_A) \to \DD_{{\rm coh}}^b (\shO_{A^{\sharp}}),
\]
where on the left we have coherent $\Dmod_A$-modules, and on the right coherent sheaves on the moduli space 
$A^{\sharp}$ of line bundle with integrable connection on $A$. This leads to generalizations of the generic vanishing 
results in \cite{GL2}.

\subsection{Generic vanishing on compact K\"ahler manifolds}

The original results on generic vanishing for $\omega_X$ were established by Green and Lazarsfeld \cite{GL1}, \cite{GL2} by means of Hodge theory, and therefore apply to compact K\"ahler manifolds. On the other hand, two well-known theorems
on projective varieties, Simpson's \cite{Simpson} result that every irreducible component of every 
$$\Sigma^k(X) = \menge{\rho \in {\rm Char}(X)}{H^k(X, \CC_{\rho}) \neq 0}$$ 
contains characters of finite order, and Hacon's result \cite{Hacon} that higher direct 
images $R^i f_* \omega_X$ are $GV$-sheaves,  have eluded on K\"ahler manifolds proofs based on similar 
methods.

This difficulty can be overcome by the use of Hodge modules. The following extension of Hacon's theorem was also motivated
by a result of Chen-Jiang \cite{CJ} for $f_* \omega_X$ on projective  varieties of maximal 
Albanese dimension.

\begin{theorem}[{\cite[Theorem A]{PPS}}]\label{thm:direct_image}
Let $f \colon X \to T$ be a holomorphic mapping from a compact K\"ahler manifold to a
compact complex torus. Then for $j \geq 0$, one has a decomposition
\[
	R^j \fl \omX \simeq \bigoplus_{k=1}^n \bigl( \qu_k \shF_k \tensor L_k \bigr),
\]
where $\shF_k$ are M-regular (hence ample) coherent sheaves with projective support 
on the compact complex tori $T_k$, $q_k \colon T \to T_k$ are surjective morphisms with connected fibers,
and $L_k \in \Pic^0(T)$ have finite order. In particular, $R^j \fl \omX$ is a
GV-sheaf on $T$.
\end{theorem}

Thus the direct images $R^j f_* \omega_X$, being direct sums of pull-backs of ample sheaves, are essentially projective objects and have positivity properties that go beyond the usual weak positivity in the sense of Viehweg (cf. \S\ref{scn:WP}).
Since they are the lowest nonzero term in filtrations on Hodge $\Dmod$-modules with integral structure, according to Example \ref{direct_image}, the result above is a consequence of the following statement generalizing Theorem \ref{hm_abelian}.
It is important to note that on K\"ahler manifolds the absence of a polarization defined over $\QQ$ leads to the need for a more general context of polarizable Hodge modules with real structure. It was already suggested in \cite{Saito-Kaehler} how to develop such a theory; for more details see \cite[Ch.B]{PPS}.

\begin{theorem}[{\cite[Theorem D]{PPS}}]\label{thm:Chen-Jiang}
Let $M = \big( (\Mmod, F); M_{\RR}\big)$ be a polarizable real Hodge
module on a complex torus $T$. Then for each $k \in \ZZ$, the coherent
$\shO_T$-module $\gr_k^F \Mmod$ decomposes as
\[
	\gr_k^F \Mmod \simeq \bigoplus_{j=1}^n 
		\bigl( \qu_j \shF_j \tensor_{\shO_T} L_j \bigr),
\]
where $q_j \colon T \to T_j$ is a surjective map with connected fibers to a 
complex torus, $\shF_j$ is an M-regular coherent sheaf on $T_j$ with projective
support, and $L_j \in \Pic^0(T)$. If $M$ admits an integral structure, then each
$L_j$ has finite order.
\end{theorem}

On the other hand, Simpson's result on points of finite order was extended to the K\"ahler setting 
by Wang \cite{Wang}, using in the process results on Hodge modules on abelian varieties from \cite{Schnell}. 

\begin{theorem}[{\cite[Theorem 1.5]{Wang}}]\label{thm:finite-order}
If a polarizable Hodge module $M$ of geometric origin on a compact complex torus $T$
admits an integral structure, then the sets
\[
	\Sigma_m^k(T, M) 
		= \menge{\rho \in {\rm Char}(T)}{\dim H^k(T, M_{\RR} \tensor_{\RR} \CC_{\rho}) \geq m}
\]
are finite unions of translates of linear subvarieties by points of finite order.
\end{theorem}

The condition of being of geometric origin essentially means that $M$ is a direct summand in the pushforward of the 
trivial complex Hodge module via a holomorphic map $f\colon X \to T$, via the Decomposition Theorem.
Wang's result was extended to arbitrary polarizable (real) Hodge modules on $T$ in \cite[Theorem E]{PPS}.

Going back to the higher direct images $R^i f_* \omega_X$, their $GV$-property can be applied to classification results for compact K\"ahler manifolds, in line with well-known statements in the projective setting. For projective varieties, part (1) below is a theorem of Chen-Hacon \cite{CH1}, answering a conjecture of Koll\'ar, while (2) is a theorem of Jiang \cite{Jiang}, providing an effective version of result of Kawamata.

\begin{theorem}[{\cite[Theorem B, Theorem 19.1]{PPS}}]\label{torus}
Let $X$ be a  compact K\"ahler manifold with $P_1 (X) = P_2(X) = 1$. Then:
\begin{enumerate}
\item If $b_1 (X) = 2 \dim X$, then $X$ is bimeromorphic to a compact complex torus.
\item More generally, the Albanese map of $X$ is surjective, with connected fibers. 
\end{enumerate}
\end{theorem}

Note that the more precise statement in Theorem \ref{thm:direct_image} leads to simplified proofs of these results even in the projective setting. Further applications to the classification of irregular manifolds with small invariants can be  found in \cite{Pareschi}.

\section{Families of varieties}

\subsection{A Hodge module construction for special families}\label{HM-construction}
One of the essential motivations for the development of Hodge modules was to extend the theory of geometric variations of Hodge structure to morphisms with arbitrary singularities. Hodge modules are thus a very useful tool for attacking problems regarding families of varieties. They provide concrete additions to this already rich area of study from two main points of view:

\begin{renumerate}
\item One can sometimes work directly with morphisms of smooth projective varieties whose singular locus is arbitrary.
\item The general fiber does not need to have strong positivity properties.
\end{renumerate}

For the statement below, we need to introduce a locus which measures the singularities 
of a surjective morphism $f \colon X \to Y$ of
smooth projective varieties. We use the following notation for the induced
morphisms between the cotangent bundles:
$$
\begin{tikzcd}
X \dar{f} & T^{\ast} Y \times_Y X \lar[swap]{p_2} \dar{p_1} \rar{\df} & T^{\ast} X \\
Y & T^{\ast} Y \lar[swap]{p}
\end{tikzcd}
$$
Inside the cotangent bundle of $Y$, we consider the set of singular cotangent vectors
\[
	S_f = p_1 \bigl( \df^{-1}(0) \bigr) \subseteq T^{\ast} Y.
\]
One can easily check that $S_f$ is the union of the zero
section and a closed conical subset of $T^{\ast} Y$ whose image via
$p$ is equal to $D_f$, the singular locus of $f$.

In the problems we consider below, starting from diverse hypotheses one always reaches the following situation:
there is a surjective morphism $f \colon X \to Y$
between two smooth projective varieties, and a line bundle $L$ on $Y$, such that if we 
consider on $X$ the line bundle
\[
	B = \omega_{X/Y} \tensor \fu L^{-1},
\]
the following condition holds: 
\begin{equation}\label{eq:sections}
	H^0 (X, B^{\otimes m}) \neq 0 \quad \text{for some $m \geq 1$.}
\end{equation}
Starting from this data, one can construct a graded module over $\shA_Y := \Sym \shT_Y$ with the
following properties.

\begin{theorem}[{\cite[Theorem 9.2]{PS2}}] \label{thm:VZ}
Assuming \eqref{eq:sections}, one can find a graded $\shA_Y$-module $\shGb$
that is coherent over $\shA_Y$ and has the following properties:
\begin{aenumerate}
\item \label{en:VZa}
As a coherent sheaf on the cotangent bundle, $\Supp \shG \subseteq S_f$.
\item \label{en:VZb}
One has $\shG_0 \simeq L \tensor \fl \OX$.
\item \label{en:VZc}
Each $\shG_k$ is torsion-free on the open subset $Y \setminus D_f$. 
\item \label{en:VZd}
There exists a regular holonomic $\Dmod_X$-module $\Mmod$ with good filtration
$F_{\bullet} \Mmod$, and an inclusion of graded $\shA_Y$-modules $\shGb \subseteq
\gr_{\bullet}^F \! \Mmod$.
\item \label{en:VZe}
The filtered $\Dmod$-module $(\Mmod, F)$ underlies a polarizable Hodge module $M$ on
$Y$ with strict support $Y$, and $F_k \Mmod = 0$ for $k < 0$.
\end{aenumerate}
\end{theorem}

The Hodge module $M$ and the graded submodule $\shGb$ are constructed by applying results from the theory of Hodge
modules to a new morphism determined by a section 
$0 \neq s \in H^0 \bigl( X, B^{\tensor m} \bigr)$. Let's at least specify $M$: such a section defines a branched covering $\pi \colon X_m \to X$ of degree $m$, unramified outside the divisor $Z(s)$. If  $m$ is chosen to be 
minimal, $X_m$ is irreducible, and we consider a resolution of
singularities $\mu \colon Z \to X_m$ that is an isomorphism over the complement of $Z(s)$. We define
$\varphi = \pi \circ \mu$ and $h = f \circ \varphi$ as in the diagram 
\begin{equation} \label{eq:geometry}
\begin{tikzcd}
Z \rar{\mu} \arrow[bend right=20]{drr}{h} \arrow[bend left=40]{rr}{\varphi} 
		& X_m \rar{\pi} & X \dar{f} \\
 & & Y
\end{tikzcd}
\end{equation}
If $d = \dim X = \dim Z$, let 
$$\shH^0 \hl \QQ_Z^H \decal{d} \in \HM(Y, d)$$ 
be the polarizable Hodge module obtained by taking the direct image of the
constant Hodge module on $Z$, provided by Theorem \ref{stability}; restricted to the smooth locus of $h$, this is just the
polarizable variation of Hodge structure on the middle cohomology of the fibers. One takes 
$M \in \HM (Y, d)$ to be the summand with strict support $Y$ in the decomposition of
$\shH^0 \hl \QQ_Z^H \decal{d}$ by strict support. Constructing $\shGb$ requires extra work, specifically using the section $s$ and the 
strictness property in Theorem \ref{stability}; see \cite[\S11]{PS2} for details.

This type of construction was introduced by Viehweg and Zuo \cite{VZ1}, \cite{VZ2}, \cite{VZ3} in the context of Higgs bundles, under the assumption that $B$ is semiample, but it becomes both simpler and more flexible through the use of Hodge modules. It is also useful for applications to use the local properties of Hodge modules in order to construct a suitable Higgs sheaf 
from the graded module $\shGb$; this is done in \cite[Theorem 9.4]{PS2}.

\subsection{Zeros of holomorphic one-forms}\label{one_forms}
The following basic result, conjectured (and proved in some important special cases) by Hacon-Kov\'acs \cite{HK} and Luo-Zhang \cite{LZ}, was obtained in \cite{PS} as an application 
of part of the construction in Theorem \ref{thm:VZ}, combined with results in \S\ref{scn:GV-projective}.


\begin{theorem}[{\cite[Theorem 2.1]{PS}}]\label{one-forms}
Let $X$ be a smooth complex projective variety, and assume that there is a linear subspace $W \subseteq H^0(X, \OmX^1)$
such that $Z(\omega)$ is empty for every nonzero $1$-form $\omega \in W$. Then 
$$\dim W \le \dim X - \kappa(X).$$
In particular, if $X$ is of general type, then every holomorphic $1$-form on $X$ has at least one zero.
\end{theorem}

It is quite standard to reduce the statement to the last assertion, i.e. to the case when $X$ is of general type.
It is also possible to express the result in the theorem in terms of the Albanese map $f\colon X \to A$ of $X$. As such, one of the key points in the proof is to have access to a Hodge module $M$ with underlying filtered $\Dmod$-module $(\Mmod, F)$
and a graded submodule $\shGb \subseteq
\gr_{\bullet}^F \! \Mmod$ on $A$ as in Theorem \ref{thm:VZ}; but in this set-up it is almost immediate to see that, at least after passing to an \'etale cover of $X$, condition $(\ref{eq:sections})$ is satisfied for any ample line bundle $L$ on $A$.

Moreover, one can show that the support of $\shGb$, seen as a coherent sheaf on $T^* A \simeq A \times H^0 (X, \Omega_X^1)$,
satisfies the property
$${\rm Supp} ~\shGb \subseteq (f \times {\rm id})(Z_f),$$
where 
$$Z_f = \{(x, \omega) ~|~ \omega (x) = 0\} \subseteq  X \times H^0 (X, \Omega_X^1).$$
Thus to prove the theorem it suffices to show that ${\rm Supp} ~\shGb$ projects surjectively onto the space of holomorphic
$1$-forms. One obtains this by using the full package of generic vanishing theorems described in \S\ref{scn:GV-projective}, applied to the Hodge module $M$.

While not directly a result about families of varieties, Theorem \ref{one-forms} has consequences to the study of families over abelian varieties 
via the following corollary, obtained by taking $W = f^* H^0 (A, \Omega_A^1)$.

\begin{corollary}[{\cite[Corollary 3.1]{PS}}]\label{cor:smooth}
If $f \colon X \to A$ is a smooth morphism from a smooth complex projective variety
onto an abelian variety, then 
$$\dim A \leq \dim X - \kappa(X).$$
\end{corollary}

This implies for example that there are no nontrivial smooth morphisms from a variety of general
type to an abelian variety; this fact was proved in \cite{VZ1} when
the base is an elliptic curve, and can be deduced from \cite{HK} in general.

\subsection{Families of varieties of general type over abelian varieties}\label{abvar}
The main application of Corollary \ref{cor:smooth} is to deduce the birational isotriviality of smooth families of varieties of general type
over abelian varieties. Note that when the fibers are canonically polarized, Kov\'acs \cite{Kovacs2} and Zhang \cite{Zhang} 
have shown that the family must be isotrivial.

\begin{theorem}[{\cite[Corollary 3.2]{PS}}]
If $f \colon X \to A$ is a smooth morphism onto an abelian variety, with 
fibers of general type, then $f$ is birationally isotrivial. 
\end{theorem}

The conclusion means that $X$ becomes birational to a product after a generically finite
base-change, i.e. that $\Var(f) = 0$ in Viehweg's terminology. It follows by combining the bound
$\kappa(X) \leq \dim F$ from \corollaryref{cor:smooth} with the solution to the $C_{n,m}^+$
conjecture for fibers of general type in \cite{Kollar2}, which in this case gives
\[
	\kappa(F) + \Var(f) \leq \kappa(X).
\]

\subsection{Families of maximal variation and Viehweg-Zuo sheaves}\label{VZ}
Arguments based on the Hodge module construction in Theorem \ref{thm:VZ} also allow, at the other end of the (variation) spectrum, to 
study families $f\colon X \to Y$ of varieties of general type (or more generally having good minimal models) with
maximal variation ${\rm Var}(f) = \dim Y$. Recall that this means that the general fiber can be birational to only countably many other fibers.

It was conjectured by Viehweg that if $Y^{\circ}$ is smooth and quasi-projective, and $f^{\circ}: X^{\circ} \rightarrow Y^{\circ}$ is a smooth family of canonically polarized varieties with maximal variation, then $Y^\circ$ must be of log-general type. More precisely, if the data is extended to compactifications $f\colon X \to Y$, and $D = Y \smallsetminus Y^{\circ}$ is a divisor, then $\omega_Y (D)$ is big.  
This conjecture in its general form was settled by Campana and P\u aun \cite{CP1}, \cite{CP2}, after a series of papers by Kov\'acs \cite{Kovacs3}, Kebekus-Kov\'acs \cite{KK1}, \cite{KK2}, \cite{KK3}, and Patakfalvi \cite{Patakfalvi} proved important special cases. One of the fundamental inputs in all of these works 
is the existence of what are called Viehweg-Zuo sheaves \cite[Theorem 1.4]{VZ2} for families of canonically polarized varieties.

The Hodge module construction in \S\ref{HM-construction} can be brought into play in order to prove the existence of Viehweg-Zuo sheaves for 
more general families.

\begin{theorem}[{\cite[Theorem B]{PS2}}]\label{thm:VZ-sheaves}
Let $f \colon X \rightarrow Y$ be an algebraic fiber space between smooth projective
varieties, such that the singular locus of $f$ is  a divisor $D_f \subseteq Y$ with simple normal crossings.
Assume that for every generically finite $\tau: \tilde Y \rightarrow Y$ with $\tilde
Y$ smooth, and for every resolution $\tilde X$ of $X\times_Y \tilde Y$, there is an
integer $m \geq 1$ such that 
$\det {\tilde f}_* \omega_{\tilde X/ \tilde Y}^{\otimes m}$ is big. Then there exists a big coherent sheaf $\shH$ on $Y$ and an integer 
$s \ge 1$, together with an inclusion
$$\shH \hookrightarrow \big(\OmY^1 (\log D_f) \big)^{\otimes s}.$$
\end{theorem}

Viehweg's $Q_{n,m}$ conjecture \cite[Remark 3.7]{Viehweg2} (see also \cite[(7.3)]{Mori}) predicts that the hypothesis of the theorem is satisfied by any algebraic fiber space $f\colon X \to Y$ 
such that $\Var (f) = \dim Y$. This conjecture was proved by Koll\'ar \cite{Kollar2} (see also \cite{Viehweg4}) when the smooth fibers of $f$ are of general type, and also by Kawamata \cite{Kawamata} when the geometric generic fiber has a good minimal model.\footnote{One of the main conjectures in birational geometry is of course that all varieties with $\kappa (X) \ge 0$ should do so.}

The point here is that the existence of an integer $m$ as in the statement of Theorem \ref{thm:VZ-sheaves} leads to the fact that 
($\ref{eq:sections}$) holds for a new family over $\tilde Y$ constructed from $f$ by means of a general strategy introduced by Viehweg, based on 
his fiber product trick and semistable reduction. Very roughly speaking, the Hodge module and graded submodule produced by Theorem \ref{thm:VZ}
can then be used in conjunction with Theorem \ref{thm:WP} on weak positivity to reach the desired conclusion.

Together with the main result of Campana-P\u{a}un \cite{CP2} on the pseudo-effectivity of quotients of powers of log-cotangent bundles, Theorem \ref{thm:VZ-sheaves} implies that Viehweg's conjecture also holds under weaker assumptions on the fibers of the family.

\begin{theorem}[{\cite[Theorem A]{PS2}}]\label{thm:hyperbolicity}
Let $f \colon X \rightarrow Y$ be an algebraic fiber space between smooth projective
varieties, and let $D \subseteq Y$ be any divisor containing the singular locus of $f$. Assume that $\Var (f) = \dim Y$. Then:
\begin{enumerate}
\item If the general fiber of $f$ is of general type, then the pair $(Y, D)$ is of
log-general type, i.e. $\omY (D)$ is big. 
\item More generally, the same conclusion holds if the geometric generic fiber of $f$ admits a good minimal model.
\end{enumerate}
\end{theorem}

\section{Hodge ideals}\label{Hodge_ideals}

\subsection{Motivation}
Let $X$ be a smooth complex variety of dimension $n$, and $D$ a reduced effective divisor on $X$. 
The left $\Dmod_X$-module
$$\shO_X (*D) = \bigcup_{k \ge 0} \shO_X (kD)$$
of functions with arbitrary poles along $D$ underlies the mixed
Hodge module $j_* \QQ_U^H [n]$, where $U = X \smallsetminus D$ and $j: U \hookrightarrow X$ is the inclusion map. It therefore comes
equipped with a Hodge filtration $F_{\bullet} \shO_X(*D)$.
Saito \cite{Saito-B} showed that this filtration is contained in the pole order filtration, namely 
$$
F_k \shO_X(*D) \subseteq \shO_X \big( (k+1) D\big) \,\,\,\,\,\, {\rm for ~all} \,\,\,\, k \ge 0,
$$
and the problem of how far these are from being different is of interest both in the study of the singularities of $D$ and in understanding 
the Hodge structure on the cohomology of the complement $H^\bullet (U, \CC)$. The inclusion 
above leads to defining for each $k \ge 0$ a coherent sheaf of ideals $I_k (D) \subseteq \shO_X$ by the formula
$$F_k \shO_X(*D) =  \shO_X \big( (k+1) D\big) \otimes I_k (D).$$
This section is devoted to a brief introduction to work with Musta\c t\u a \cite{MP}, in which we study and apply the theory of these ideals using birational geometry methods. This involves redefining them by means of log-resolutions. 

\subsection{Alternative definition and local properties} 
We consider a log-resolution $f \colon Y \to X$ of the pair $(X, D)$, and denote $E = \tilde D + F$, where $F$ is the 
reduced exceptional divisor. Because we are dealing with pushforwards, it is more convenient to work in the setting of right 
$\Dmod$-modules; in other words we will look at a filtration on $\omega_X (*D)$.

The starting point is the existence of a filtered complex of right $f^{-1} \Dmod_X$-modules 
$$0 \longrightarrow f^* \Dmod_X \longrightarrow \Omega_Y^1(\log E) \otimes_{\shO_Y} f^* \Dmod_X \longrightarrow \cdots $$
$$\cdots \longrightarrow \Omega_Y^{n-1}(\log E) \otimes_{\shO_Y} f^* \Dmod_X \longrightarrow \omega_Y(E) \otimes_{\shO_Y} 
f^*\Dmod_X \longrightarrow 0$$
which is exact except at the rightmost term, where the cohomology is $\omega_Y(*E) \otimes_{\Dmod_Y} \Dmod_{Y\to X}$. 
Denoting it by $A^\bullet$, it has a filtration with subcomplexes $F_{k-n}A^\bullet$, for $k \ge 0$,
given by 
$$0 \rightarrow f^* F_{k-n} \Dmod_X \rightarrow \Omega_Y^1(\log E) \otimes_{\shO_Y} f^* F_{k-n+1} \Dmod_X \rightarrow \cdots \to  \omega_Y(E) \otimes_{\shO_Y} f^* F_k \Dmod_X\rightarrow 0.$$

\begin{definition}
The filtration on $\omega_X(*D)$, and in particular the \emph{$k$-th Hodge ideal} $I_k(D)$ associated to $D$, are defined by the formula
$$F_k \omega_X(*D) : = \omega_X \big( (k +1)D \big) \otimes I_k(D) = {\rm Im} \left[R^0 f_* F_{k-n}A^\bullet \longrightarrow R^0 f_* A^\bullet \right],$$
after proving that this image is contained in $\omega_X \big( (k +1)D \big)$. It is shown in \cite{MP} that this definition is independent of the 
choice of log resolution, and that it indeed coincides with Saito's Hodge filtration, as in the previous paragraph. 
\end{definition}

\begin{example}
(i) If $D$ is smooth, then $I_k (D) = \shO_X$ for all $k \ge 0$.

\noindent
(ii) More generally, when $D$ is a simple normal crossings divisor, one has  
$$\omega_X \big( (k+1)D\big) \otimes I_k (D) = \omega_X (D) \cdot F_k \Dmod_X \,\,\,\,\,\, {\rm  for~ all} \,\,\,\, k \ge 0.$$

\noindent
(iii) For $k=0$ and $D$ arbitrary, we have 
$$I_0 (D) = \mathcal{I} \big(X, (1-\epsilon)D \big),$$
the multiplier ideal associated to the $\QQ$-divisor $(1-\epsilon)D$ on $X$, for any $0 < \epsilon \ll 1$. In particular, $I_0 (D) = \shO_X$ if and only if $(X, D)$ is log-canonical; see \cite[9.3.9]{Lazarsfeld}.
\end{example}

This last example explains one of the main motivations for studying the higher ideals $I_k (D)$; taken together, they provide 
a generalization of the notion of multiplier ideal, and a finer invariant of singularities. It turns out that there is a sequence  of inclusions
$$\cdots I_k (D) \subseteq \cdots \subseteq I_1 (D) \subseteq I_0 (D).$$
(Note that the fact that $F_k \omega_X (*D)$ is a filtration in the sense of $\Dmod_X$-modules yields reverse inclusions
$I_k (D) \cdot \shO_X(-D) \subseteq I_{k+1}(D)$.)

\begin{definition}
The pair $(X, D)$ is \emph{$k$-log-canonical} if $I_k (D) = \shO_X$.
Thus $0$-log-canonical is equivalent to the usual notion of a log-canonical pair.  The general notion refines the property of having rational 
singularities:
\end{definition}

\begin{theorem}[{\cite[Theorem B]{MP}}]\label{rational}
For every $k \ge 1$ we have an inclusion
$$I_k(D) \subseteq {\rm adj}(D),$$
where ${\rm adj}(D)$ is the adjoint ideal of the divisor $D$; see \cite[\S9.3.E]{Lazarsfeld}.
Hence if $(X, D)$ is $k$-log canonical for some $k \ge 1$, then $D$ is normal with rational singularities.
\end{theorem}

Recall from the example above that if $D$ is smooth, all Hodge ideals are trivial. It turns out that the two conditions are equivalent; in fact, more precisely, 
any level of log-canonicity beyond $(n-1)/2$ implies smoothness.

\begin{theorem}[{\cite[Theorem A]{MP}}]\label{smoothness}
The following are equivalent:

\noindent
(i) $D$ is smooth. 

\noindent
(ii) the Hodge filtration and pole order filtration on $\shO_X(*D)$ coincide.

\noindent
(iii) $I_k (D) = \shO_X$ for all $k \ge 0$.

\noindent
(iv) $I_k (D) = \shO_X$ for some $k \ge \frac{n-1}{2}$.
\end{theorem}

The statements above are consequences of general results regarding the order of vanishing of $I_k(D)$ along exceptional divisors on 
birational models over $X$, and along closed subsets in $X$.  This in particular applies to provide criteria for the nontriviality of the Hodge ideals in terms of the multiplicity of $D$ at singular points; one such criterion is the following:

\begin{theorem}[{cf. \cite[Theorem E]{MP}}]\label{criterion_nontriviality}
Let $D$ be a reduced effective divisor on a smooth variety $X$, and let $x\in X$. If $m={\rm mult}_x(D)$, then for every $k$ we have
$$I_k(D)\subseteq \frak{m}_x^q,\quad\text{where}\text\quad q=\min\{m-1,(k+1)m-n\},$$
where $\frak{m}_x^q=\shO_X$ if $q\leq 0$.
\end{theorem}

As an example, this says that for $q \le m-1$ one has 
$$m \ge \frac{n+q}{2} \implies I_1(D) \subseteq \frak{m}_x^q.$$

\subsection{Examples}
Here are  some concrete calculations or statements that can be made at the moment. It is useful to start by introducing another related class of ideals.
By analogy with the simple normal crossings case, one can define for each $k \ge 0$ auxiliary ideal sheaves $J_k (D)$ by the formula
$$\shO_X \big((k+1)D\big) \otimes J_k (D) =  F_k \Dmod_X \cdot \big( \shO_X (D) \otimes I_0 (D)\big) = F_k \Dmod_X \cdot F_0 \shO_X (*D) .$$
Since $F_k \shO_X(*D)$ is a filtration in the sense of $\Dmod_X$-modules, it is clear that for each $k \ge 0$ there is an inclusion
$$J_k(D) \subseteq I_k(D).$$
In dimension at least three this inclusion is usually strict, and the Hodge ideals are much more difficult to compute than $J_k (D)$. For instance, for $k =1$, 
using the log resolution notation above, one can verify the existence of a short exact sequence 
$$0 \longrightarrow \omega_X (2D) \otimes J_1(D) \longrightarrow \omega_X (2D) \otimes I_1(D)  \longrightarrow R^1 f_* \Omega_Y^{n-1} (\log E) 
\longrightarrow 0.$$

\begin{example}[{\bf Surfaces}]
When $\dim X = 2$ however, and more generally outside of a closed subset of codimension $\ge 3$ in any dimension, it is shown in \cite{MP} that $I_k (D) = J_k (D)$ 
for all $k \ge 0$. Here are some concrete calculations:

\medskip

\noindent
$\bullet$~~ If $D = (xy = 0) \subset \CC^2$ is a node, then $I_k (D) = (x, y)^k$ for all $k \ge 0$.

\noindent
$\bullet$~~ If $D = (x^2 + y^3 = 0) \subset \CC^2$ is a cusp, then $I_0(D) = (x, y)$, $I_1 (D) = (x^2, xy, y^3)$, and $I_2 (D) = (x^3,x^2y^2,xy^3,y^5,y^4-3x^2y)$.

\noindent
$\bullet$~~ If $D = (xy(x+y)= 0) \subset \CC^2$ is a triple point, then $I_0(D) = (x, y)$, while $I_1 (D) = (x,y)^3$.

Note for instance how $I_1$ distinguishes between singularities for which $I_0$ is the same. This is one of the ways in which Hodge ideals 
become important in applications. For details on all of the examples that come next, see \cite[\S20]{MP}.
\end{example}

\begin{example}[{\bf Ordinary singularities}]
Let $x \in D$ be a point of multiplicity $m \ge 2$, with the property that the projectivized tangent cone of $D$ at $x$ is smooth; these are sometimes called 
ordinary singularities. For instance, $D$ could be the cone over a smooth hypersurface of degree $m$ in $\AAA^{n-1}$. One can show:
\medskip

\noindent
$\bullet$~~ $I_k (D)_x = \shO_{X,x} \iff k \le \left[\frac{n}{m}\right] - 1$.

\noindent
$\bullet$~~If $\frac{n}{2}\le m\le n-1$, then $I_1 (D)_x = \frak{m}_x^{2m- n}$ .

\noindent
$\bullet$ ~~If $m \ge n$, then
$$\shO_X(-D)\cdot \frak{m}_x^{m-n-1}+\frak{m}_x^{2m-n}\subseteq I_1(D)\subseteq \shO_X(-D)\cdot \frak{m}_x^{m-n-2}+\frak{m}_x^{2m-n-1},$$ 
with $\dim_{\CC}I_1(D)/(\shO_X(-D)\cdot \frak{m}_x^{m-n-1}+\frak{m}_x^{2m-n})=m{{m-2}\choose{n-2}}$. 
\end{example}

\begin{example}[{\bf Determinantal varieties and theta divisors}]
Regarding Theorem \ref{rational}, we see already from the example above that there can be divisors $D$ with rational singularities such that $I_k (D) \neq \shO_X$ for all $k \ge 1$. Here are two more celebrated examples where things can go both ways:

\medskip

\noindent
$\bullet$~~If $(A, \Theta)$ is an irreducible principally polarized abelian variety, a result of Ein-Lazarsfeld \cite{EL} says that $\Theta$ is normal, with 
rational singularities, i.e. ${\rm adj}(\Theta) = \shO_A$. On the other hand, it can happen that $I_1 (\Theta) \neq \shO_A$. For example, when $(A, \Theta)$
is the intermediate Jacobian of a smooth cubic threefold, it follows from the example above that $I_1(\Theta) = \frak{m}_0$; indeed, 
the origin $0 \in A$ is the unique point singular point of $\Theta$, and is an ordinary singularity of multiplicity $3$. 

\noindent
$\bullet$~~If $D\subset X = \AAA^{n^2}$ is the generic determinantal variety of codimension $1$, one can show that $I_1 (D) = \shO_X$. This shows that 
$D$ has rational singularities, which is of course a well-known result. Moreover, it turns out that  $D$ is $1$-log-canonical but not more, 
i.e. $I_2 (D) \neq \shO_X$.
\end{example}

\begin{example}[{\bf Roots of Bernstein-Sato polynomial and diagonal hypersurfaces}]
Suppose that $f$ is a local equation of $D$. Recall that the Bernstein-Sato polynomial of $D$ at $x$ is the monic polynomial $b_{f,x}\in\CC[s]$
of smallest degree such that around $x$ there is a relation
$$b_{f,x}(s)f^s=P(s)\bullet f^{s+1}$$
for some nonzero $P\in \Dmod_X[s]$. It is known that $(s+1)$ divides $b_{f,x}$ and all roots of $b_{f,x}(s)$ are negative rational numbers. One defines
$\alpha_{f,x}$ to be $-\lambda$, where $\lambda$ is the largest root of $b_{f,x}(s)/(s+1)$. Saito showed in \cite[Theorem~0.11]{Saito-B} that 
$$I_k(D)_x=\shO_{X,x}\quad\text{for all}\quad k\leq \alpha_{f,x}-1.$$
When $D$ is the divisor in $\AAA^n$ defined by $f=\sum_{i=1}^nx_i^{a_i}$, with $a_i\geq 2$, it is known that 
$\alpha_{f,0}=\sum_{i=1}^n\frac{1}{a_i}$, and so it follows that
$$I_k(D)=\shO_X\quad\text{for all}\quad k\leq-1+\sum_{i=1}^n\frac{1}{a_i}.$$
\end{example}

\subsection{Vanishing and applications}
As it is well understood from the study of multiplier ideals, criteria like Theorem \ref{criterion_nontriviality} are most useful for applications when applied in combination with vanishing theorems. There is indeed such a general vanishing theorem for Hodge ideals. For $k = 0$, it is precisely the well-known Nadel Vanishing; for $k \ge 1$ it requires more careful hypotheses. 

\begin{theorem}\label{higher_vanishing}
Let $X$ be a smooth projective variety of dimension $n$, $D$ a reduced effective divisor, and $L$ a line bundle on $X$. Then, for each $k \ge 1$, assuming that the pair $(X,D)$ is $(k-1)$-log-canonical we have:

\medskip

\noindent
(i) If $L$ is a line bundle such that $L (p D)$ is ample for all $0 \le p \le k$,\footnote{If $k \ge \frac{n+1}{2}$, then $D$ is in fact smooth by 
Theorem \ref{smoothness}, and so $I_k (D) = \shO_X$.
 In this case, if $L$ is a line bundle such that $L (k D)$ is ample, then 
$$H^i \big(X, \omega_X ( (k+1)D) \otimes L  \big) = 0 \,\,\,\,\,\, {\rm for ~all}\,\,\,\, i >0.$$} then 
$$H^i \big(X, \omega_X ( (k+1)D) \otimes L \otimes I_k (D) \big) = 0$$
for all $i \ge 2$. Moreover, 
$$H^1 \big(X, \omega_X ((k+1)D) \otimes L \otimes I_k (D) \big) =  0$$
holds if $H^j \big(X, \Omega_X^{n-j} \otimes L ((k - j +1)D )\big) = 0$ for all $1 \le j \le k$.

\medskip

\noindent
(ii) If $D$ is ample, then (1) and (2) also hold with $L = \shO_X$.
\end{theorem}

The main ingredient in the proof is the Kodaira-Saito vanishing theorem for Hodge modules, Theorem \ref{saito_vanishing}.
In important classes of examples, the extra hypotheses in Theorem \ref{higher_vanishing} are either automatically satisfied, or can be discarded, which
greatly increases the range of applications. For instance, in all of the examples below the borderline Nakano-type hypothesis in (i) holds automatically:

\medskip

\noindent
$\bullet$~~\emph{Toric varieties}: due to the Bott-Danilov-Steenbrink vanishing theorem.

\noindent
$\bullet$~~\emph{Projective space} $\PP^n$: in addition to the above, one can also discard the assumption that $D$ be $(k-1)$-log canonical, essentially because of the
existence of a Koszul resolution for all $\Omega_{\PP^n}^j$.

\noindent
$\bullet$~~\emph{Abelian varieties}: due to the fact that all $\Omega_X^j$ are trivial; for the same reason, one can also discard the hypothesis that 
$D$ be $(k-1)$-log canonical.

\medskip

As mentioned above, when combined with nontriviality criteria like Theorem \ref{criterion_nontriviality}, these stronger versions of the vanishing theorem have interesting consequences regarding the behavior of isolated singular points on hypersurfaces in $\PP^n$, or on theta divisors on ppav's. Here is a sampling of applications from \cite{MP}; for further results and context see Ch.H in \emph{loc. cit.}

\begin{theorem}
Let $D$ be a reduced hypersurface of degree $d$ in $\PP^n$, with $n \ge 3$, and denote by $S_m$ the set of isolated singular points on $D$ of multiplicity $m\ge 2$. Then $S_m$ imposes independent conditions on 
hypersurfaces of degree at least $([\frac{n}{m}]+ 1)d - n -1$.
\end{theorem}

This turns out to be an improvement of the bounds found in the literature for most $n\ge 5$ and $m\ge 3$.

A well-known result of Koll\'ar \cite[Theorem 17.3]{Kollar4} states that if $(A, \Theta)$ is a principally polarized abelian variety (ppav) of dimension $g$, then the pair $(A, \Theta)$ is log-canonical, and so ${\rm mult}_x (\Theta) \le g$ for any $ x \in \Theta$. It is well known however that when $\Theta$ is 
irreducible one can do better; see e.g. \cite{EL} and the references therein. It is a folklore conjecture that the bound should actually be $\frac{g+1}{2}$. 

\begin{theorem}\label{intro_theta_isolated}
Let $(A, \Theta)$ be ppav of dimension $g$ such that $\Theta$ has isolated singularities. 
Then:

\medskip

\noindent
(i)  For every $x \in \Theta$ we have
${\rm mult}_x (\Theta) \le \frac{g + 1}{2}$, and also ${\rm mult}_x (\Theta) \le \epsilon(\Theta) + 2$, where $\epsilon (\Theta)$ is the Seshadri constant of the principal polarization.

\noindent
(ii) Moreover, there can be at most one point $x \in \Theta$ such that ${\rm mult}_x (\Theta) = \frac{g + 1}{2}$.
\end{theorem}

This is obtained using the ideal $I_1(\Theta)$, except for the bound involving the Seshadri constant, which uses all $I_k(\Theta)$.
Results of a similar flavor have also been obtained recently by Codogni-Grushevsky-Sernesi.

\section{Applications of Hodge modules in other areas}\label{other_applications}

While my lack of expertise does not allow me to go into much detail, I will at least try to indicate some other areas in which constructions based on Hodge modules have made a substantial impact in recent years. The list is by no means supposed to be 
exhaustive; I apologize for any omissions.  Many more references regarding the topics below, as well as others, can be found in \cite[\S3.7]{Saito-YPG}.

\noindent
{\bf Hodge theory.}
By its very nature, M. Saito's study of Hodge modules \cite{Saito-MHP}, \cite{Saito-MHM} has initially provided fundamental applications to core Hodge theory, like for instance the existence of mixed Hodge structures on the cohomology of complex varieties or intersection cohomology, or 
the enhancement of the Decomposition Theorem involving Hodge filtrations. More recently, Hodge modules have found new applications internal to  Hodge theory, 
especially to the study of (admissible) normal functions. In \cite{BFNP}, the authors used mixed Hodge modules to show
that the Hodge conjecture is equivalent to the existence of singularities in certain admissible normal functions.
The article \cite{Schnell4} contains the construction of complex analytic N\'eron models that extend families of intermediate Jacobians associated to complex variations of Hodge structure of odd weight over Zariski open subsets of complex manifolds; Hodge modules 
provide a functorial construction, with no extra requirements at the boundary. 

\noindent
{\bf Singularities, characteristic classes.} 
It also goes without saying that Hodge modules have made an impact in the study of singularities via Hodge theoretic methods, not 
far from the spirit of Ch.\ref{Hodge_ideals} above. Besides \cite{Saito-B} invoked in the text, I will highlight 
the papers \cite{Saito-Spectrum}, which proves a conjecture of Steenbrink on the spectrum of certain isolated singularities, 
\cite{BS}, which makes a connection between multiplier ideals and the $V$-filtration, and \cite{DMS}, which extends the latter to arbitrary subvarieties; see also the references therein. Hodge modules have also proved 
useful in the study of characteristic classes of singular varieties; see the survey \cite{MS}.

\noindent
{\bf Donaldson-Thomas theory.}
Mixed Hodge modules have recently been used successfully in Donaldson-Thomas theory. For instance, the cohomological $DT$ invariants 
(the hypercohomology of a certain perverse sheaves of vanishing cycles on  moduli spaces of simple sheaves on Calabi-Yau threefolds) can be 
endowed with mixed Hodge structures via Hodge module theory. They are also used in the categorification of Donaldson-Thomas invariants. 
I will highlight \cite{KL} and  \cite{BBDJS}, and refer to \cite{Szendroi} for a survey and for an extensive list of references.

\noindent
{\bf Representation theory.}
Hodge modules have become an important tool in geometric representation theory. They have been used in the study of unitary representations of reductive Lie groups \cite{SV}. A categorical ${\rm sl}_2$-action on a special category of sheaves on Grassmannians is related to the calculation in \cite{CDK} of the associated graded with respect to the Hodge and weight filtration 
for a representation-theoretically relevant example of localization along a hypersurface. See also \cite{AK} a relationship with Koszul duality.

\noindent
{\bf Automorphic forms.}
In \cite{Suh}, J. Suh has been able to apply a vanishing theorem for the associated graded 
quotients of de Rham complexes of Deligne canonical extensions, proved using tools from Hodge module theory (and the main ingredient in his proof of Theorem \ref{KV})  in order to obtain  vanishing theorems for the cohomology of Shimura varieties. The use of Hodge modules has allowed for treating the coherent cohomology groups of automorphic bundles on arbitrary Shimura varieties, with general weights, extending previous results of Faltings, Li-Schwermer, and Lan-Suh.

\bibliographystyle{amsalpha}
\bibliography{bibliography}

\end{document}